\documentclass{amsart}
\usepackage[utf8]{inputenc}
\usepackage{microtype}
\usepackage{amssymb,mathtools,mathrsfs}
\usepackage[hidelinks]{hyperref}
\usepackage{bookmark}
\hypersetup{pdftitle={Finite rank and compactness of finite sums of Toeplitz products on the harmonic Bergman space},pdfauthor={Puyu Cui and Zengyong Huang}}
\allowdisplaybreaks[1]
\theoremstyle{plain}
\newtheorem{theorem}{Theorem}[section]
\newtheorem{lemma}[theorem]{Lemma}
\newtheorem{proposition}[theorem]{Proposition}
\newtheorem{corollary}[theorem]{Corollary}
\theoremstyle{remark}
\newtheorem{remark}[theorem]{Remark}
\theoremstyle{definition}

\newcommand{\D}{\mathbb D}
\newcommand{\C}{\mathbb C}
\newcommand{\HH}{\mathbb H}
\newcommand{\A}{A^2(\D)}
\newcommand{\Azer}{zA^2(\D)}
\newcommand{\bh}{b^2(\D)}
\newcommand{\hinfty}{h^\infty(\D)}

\DeclareMathOperator{\rank}{rank}
\newcommand{\TT}{\mathsf T}
\newcommand{\Loc}{C^\infty_{\mathrm{loc}}}
\numberwithin{equation}{section}
\begin{document}

\title[Finite sums of Toeplitz products]{Finite rank and compactness of finite sums of Toeplitz products on the harmonic Bergman space}
\author[P. Cui]{Puyu Cui}
\address{School of Mathematics, Liaoning Normal University, Dalian 116029, China}
\email{cuipuyu1234@163.com}
\author[Z. Huang]{Zengyong Huang*}
\address{School of Mathematics, Liaoning Normal University, Dalian 116029, China}
\email{huangzengyong123@163.com}
\thanks{* Corresponding author.}
\date{}
\begin{abstract}
We completely characterize finite rank and compactness for finite sums
of products of Toeplitz operators on the harmonic Bergman space
of the unit disk, assuming boundedness only of the full harmonic symbols.
Reflecting the appropriate symbols reduces sums of Toeplitz products
and small Hankel products to Toeplitz product sums. For sums of mixed
Toeplitz and small Hankel products in either order, right multiplication
by a fixed unitary reflection gives the corresponding reduction.
These reductions preserve rank and compactness.

Finite rank is characterized by symbol identities and a cubic divisibility condition on the coefficient polynomials obtained after choosing bases over one-variable rational function fields.  We also obtain
an exact rank formula. Compactness is characterized by boundary
vanishing of the diagonal symbol sum and uniform decay of the centered
symbol tensor on hyperbolic disks of fixed radius as their centers
approach the boundary. For a sum of $N\ge1$ products, an equivalent
compactness criterion uses weighted derivatives of the symbol components
up to order $N$. When $N$ is even, the conditions of order $N$ cannot
be omitted. Mixed commutators of finite rank vanish. We characterize
small Hankel commutators of nonzero finite rank in terms of a bounded
analytic function and its bounded reciprocal, both with rational
primitives.
\end{abstract}
\maketitle
\raggedbottom

\section{Introduction}\label{sec:intro}

We study finite sums of Toeplitz products with bounded harmonic
symbols on the harmonic Bergman space. 
For a function $f$ on $\D$, write $(Uf)(z)=f(\bar z)$ for reflection
and $\overline f(z)=\overline{f(z)}$ for pointwise conjugation.
The small Hankel operator with symbol $\alpha$ is $h_{\alpha}=T_{\alpha}U$ \cite{GZ}.
There are two groups of products:
a sum of $TT$ and $hh$ terms is a Toeplitz product sum after reflection
of the appropriate symbols; a sum of $Th$ and $hT$ terms becomes one
after right multiplication by the same unitary $U$.
Both reductions preserve rank and compactness
(Remark~\ref{rem:product-types}). For bounded harmonic symbols, our
criteria apply to these groups separately. They do not classify an
arbitrary mixture $D_0+D_1U$, where $D_0,D_1$ are Toeplitz product sums.

On the analytic Bergman space, Luecking's theorem \cite{Luecking}
rules out nonzero finite-rank Toeplitz operators with bounded function
symbols. Guo, Sun, and Zheng \cite{GSZ} established finite-rank rigidity
for commutators and semicommutators with bounded harmonic symbols, and
for sums of products of the form $H_g^*H_f$ of big Hankel operators.
Choe, Koo, and Lee \cite{CKL} characterized when a finite sum of
Toeplitz products with bounded harmonic symbols on the analytic
Bergman space vanishes, has finite rank, or is compact.
For mixed commutativity on the analytic space, Wang, Zhang, and Zhao
\cite{WZZ} treated Toeplitz symbols of the form $\bar p+\varphi$, with
$p$ polynomial and $\varphi$ bounded analytic, and bounded analytic
small Hankel symbols. Axler and Zheng \cite{AZ} proved that a finite
sum of finite Toeplitz products with bounded symbols is compact
precisely when its Berezin transform vanishes at the boundary.

Related results on the Hardy, Dirichlet and Fock spaces can be found in
\cite{LeeHardy,LeeDirFinite,LeeDirCompact,BI}.

The harmonic Bergman space contains both analytic and coanalytic
functions. For example, the expression $T_z^2-T_{z^2}T_1$ vanishes
on the analytic Bergman space. On the harmonic Bergman space it has
rank one and sends $\bar z$ to $-z/6$, as in the monomial family of
Dong and Zhou \cite{DZ} (see \eqref{eq:two-orders-rankr-def} below).
Guo and Zheng \cite{GZ} used the analytic and coanalytic block
decomposition to study the Toeplitz and Hankel algebras on the harmonic
Bergman space and described the corresponding quotient algebras for
symbols continuous on the closed disk. Choe, Koo, and Na \cite{CKNa}
obtained finite-rank rigidity for Toeplitz products on harmonic Bergman
spaces of real balls under boundary regularity assumptions. 

Guan and Lu \cite{GL} studied products and commutativity of Toeplitz
and small Hankel operators with radial and quasihomogeneous symbols.
Dong and Zhou \cite{DZ} obtained finite-rank classifications and rank
formulas for classes of quasihomogeneous commutators and generalized
semicommutators. Quasihomogeneity constrains the angular dependence but
allows radial factors outside the harmonic class; these hypotheses
and bounded harmonicity are therefore not comparable in general.
For bounded harmonic symbols, Zhang and Zhao \cite{ZZ} obtained
necessary conditions and sufficient conditions for mixed commutativity,
with equivalent characterizations when the Toeplitz symbol is analytic
or when the harmonic Hankel symbol has real Taylor coefficients.
These are specific predecessors for the mixed commutator considered
in Section~\ref{sec:examples}.

This paper continues the study of bounded harmonic symbols in
Cui, Huang, Lu, and Zu \cite{CHLZ}, which classifies finite-rank
Toeplitz commutators and characterizes their compactness by weighted
derivative conditions. We use its weighted kernel identities and affine
rescaling method, with specific references where they enter the proofs.
The additional steps here concern cancellations between arbitrarily many
product terms. Lemma~\ref{lem:quotient} computes the ranks of the
combined kernels from their coefficients over rational function fields.
Proposition~\ref{prop:halfplane-zero} characterizes general zero sums on
the half-plane by an additive identity for the full symbol tensor.
Applied to boundary limits, it gives a compactness criterion involving
derivatives of each symbol component up to order $N$.
Remark~\ref{rem:compact-even-order} shows that this bound is sharp for
even $N$.

Compactness in harmonic spaces requires particular care. Choe and Nam
\cite{CN} studied finite Toeplitz product sums on harmonic Bergman
half-spaces and explained why the Berezin transform alone cannot give
a general compactness test: it annihilates every commutator of Toeplitz operators. Liu and Lin \cite{LiuLin} constructed a weakly localized operator algebra containing bounded-symbol Toeplitz operators and
gave a compactness criterion using normalized kernels at nearby
points. This local criterion is an operator-theoretic input to our
boundary analysis.

Let $dA=\pi^{-1}dx\,dy$ on the unit disk $\D$, and let $Q$ be the
orthogonal projection of $L^2(\D,dA)$ onto the complex harmonic Bergman
space $\bh$. Write $\A$ for the analytic Bergman space and $\Azer$
for its subspace of functions vanishing at the origin. Let $P$, $P_0$,
$\bar P$, and $\bar P_0$ be the orthogonal projections from
$L^2(\D,dA)$ onto $\A$, $\Azer$, $\overline{\A}$, and
$\overline{\Azer}$, respectively.  With $E_0u=\langle u,1\rangle1$, we have
\[
 P=P_0+E_0,\qquad \bar P=\bar P_0+E_0,\qquad Q=P+\bar P_0.
\]
The analytic Bergman kernel is $K_w(z)=\frac{1}{(1-z\bar w)^2}$, $z, w \in \D$.

For  $\alpha \in L^\infty (\D, dA)$, let
\[T_{\alpha}=QM_{\alpha}|_{\bh}, \TT_{\alpha}=PM_{\alpha}|_{\A}\] denote the harmonic and
analytic Toeplitz operators, respectively. Throughout this paper, we denote by $M_{\alpha}$ the multiplication operator with symbol $\alpha$. Let $h^\infty(\D)$ be
the space of bounded complex-valued harmonic functions.

For an integer $N\ge1$, consider
\begin{equation}\label{eq:S}
 D=\sum_{j=1}^{N}T_{\varphi_j}T_{\psi_j},
 \qquad \varphi_j,\psi_j\in h^\infty(\D).
\end{equation}

\subsection{Finite-rank characterization}\label{subsec:finite-main}

For the symbols of $D$ in \eqref{eq:S}, write the unique centered
decompositions
\begin{equation}\label{eq:centered}
\begin{aligned}
\varphi_j(z)&=\varphi_j(0)+f_j(z)+g_j(\bar z),\\
\psi_j(z)&=\psi_j(0)+p_j(z)+q_j(\bar z),
\qquad 1\le j\le N.
\end{aligned}
\end{equation}
The four analytic components vanish at zero. The
boundedness of the full harmonic symbols implies that their analytic
components belong to $H^2(\D)$ \cite[Section~3]{GSZ}.

For analytic $f,h$, let $F_f,F_h,F_{fh}$ be the primitives of
$f,h,fh$, respectively, with value zero at the origin. Define
\begin{equation}\label{eq:Rdef}
\begin{split}
\mathscr R_{f,h}(z,w)
={}&(z-w)f(z)h(w)
-f(z)\bigl(F_h(z)-F_h(w)\bigr)\\
&-h(w)\bigl(F_f(z)-F_f(w)\bigr)
+F_{fh}(z)-F_{fh}(w).
\end{split}
\end{equation}
 Moreover,
\begin{equation}\label{eq:Rintegral}
\mathscr R_{f,h}(z,w)
=\int_w^z\bigl(f(t)-f(z)\bigr)\bigl(h(t)-h(w)\bigr)\,dt,
\end{equation}
so $\mathscr R_{f,h}(z,w)/(z-w)^3$ extends holomorphically
across the diagonal.

We next specify a finite test for any finite separated expression
$\mathcal{R}(z,w)=\sum_{\ell=1}^N \alpha_\ell(z)\gamma_\ell(w)$.
In the meromorphic function field, choose bases $\xi_1,\ldots,\xi_r$
and $\eta_1,\ldots,\eta_m$ of the spans of its left and right factors
over $\C(z)$ and $\C(w)$, respectively. Write
\begin{equation}\label{eq:test-expansion}
\begin{aligned}
\mathcal{R}(z,w)&=\sum_{i=1}^{r}\sum_{j=1}^{m}
\mathcal{R}_{ij}(z,w)\xi_i(z)\eta_j(w),\\
\mathcal{R}_{ij}(z,w)&=\frac{P_{ij}(z,w)}{A(z)B(w)}.
\end{aligned}
\end{equation}
Here $A,B$ are common nonzero one-variable polynomial denominators,
and $P_{ij}\in\C[z,w]$. Rationality refers to the coefficients
$\mathcal{R}_{ij}$; the basis functions and the analytic components in
\eqref{eq:centered} need not be rational. The required test is
\begin{equation*}\tag{RD3}\label{eq:rd3}
(z-w)^3\mid P_{ij}(z,w)
\qquad\text{for every }i,j.
\end{equation*}
Equivalently, for every coefficient polynomial,
\[
P_{ij}(z,z)=0,\qquad
(\partial_wP_{ij})(z,z)=0,\qquad
(\partial_w^2P_{ij})(z,z)=0.
\]
The zero expression satisfies this condition. Lemma~\ref{lem:quotient}
shows that the outcome does not depend on the finite separated
expression, the rational bases, or the denominators chosen.
Once the rational bases are known, \eqref{eq:rd3} is a finite polynomial test.
Finding them for a given analytic family may require a separate argument.

If \eqref{eq:rd3} holds, expand
\[
\frac{P_{ij}(z,w)}{(z-w)^3}
=\sum_{s,t\ge0}c_{is,jt}z^s w^t,
\]
and let $M_\mathcal{R}=(c_{is,jt})$ be the resulting
finite coefficient matrix, with zero rows and columns omitted.
Set $\rank M_0=0$.

\begin{theorem}\label{thm:general}
Let $D$ be given by \eqref{eq:S}, with centered components as in \eqref{eq:centered}.
Then $D$ has finite rank if and only if the following three conditions hold:
\begin{align}
 \textnormal{(i)}\quad &\sum_{j=1}^{N}\varphi_j(z)\psi_j(z)=0
       &&(z\in\D),\label{eq:full-symbol-sum}\\
 \textnormal{(ii)}\quad &\sum_{j=1}^{N} f_j(z)q_j(w)=0
       &&(z,w\in\D),\label{eq:mixed-identities}\\
 \textnormal{(iii)}\quad &\mathcal{R}_+=\textstyle\sum_{j=1}^{N}\mathscr R_{f_j,p_j},\quad
        \mathcal{R}_-=\textstyle\sum_{j=1}^{N}\mathscr R_{g_j,q_j}
       &&\text{both satisfy \eqref{eq:rd3}.}\label{eq:two-numerators}
\end{align}
In that case $D1=0$, $D(\Azer)\subset\overline{\Azer}$ and
$D(\overline{\Azer})\subset\Azer$. Identifying $\overline{\Azer}$
with $\Azer$ by $U$, we have
\begin{align}
 \rank D&=\rank M_{\mathcal{R}_+}+\rank M_{\mathcal{R}_-},\label{eq:rank-general}\\
 D&\simeq0_{\C}\oplus\begin{pmatrix}0&C_+\\ C_-&0\end{pmatrix},
       \qquad \rank C_\pm=\rank M_{\mathcal{R}_\pm}.\label{eq:general-diagonal-blocks}
\end{align}
Together with the zero operator, the monomial examples of Dong and Zhou
\cite{DZ} realize every nonnegative integer as the rank of a two-term sum.
\end{theorem}

\subsection{Compactness characterizations}\label{subsec:compact-main}

For the Toeplitz sum $D$ in \eqref{eq:S}, define
\begin{equation}\label{eq:Phi}
 \Phi(z,w)=\sum_{j=1}^{N}\varphi_j(z)\psi_j(w).
\end{equation}
Set
\[
 \beta(z,w)=\operatorname{arctanh}
    \left|\frac{z-w}{1-\bar wz}\right|,
 \qquad D_\beta(a,R)=\{z:\beta(z,a)<R\},
\]
and, for a fixed $R>0$, define
\begin{equation}\label{eq:omega}
 \omega_R(a)=\sup_{z,w\in D_\beta(a,R)}
 \left|\Phi(z,w)-\Phi(z,a)-\Phi(a,w)+\Phi(a,a)\right|.
\end{equation}

Choose analytic decompositions of the bounded harmonic symbols in \eqref{eq:S}:
\begin{equation}\label{eq:components}
 \varphi_j=F_j+\overline{G_j},\qquad
 \psi_j=H_j+\overline{K_j}.
\end{equation}
These capital-letter components are used for the compactness criteria;
\eqref{eq:centered} uses the centered reflection decomposition instead.

\begin{theorem}\label{thm:compact}
For the Toeplitz sum $D$ in \eqref{eq:S}, use
\eqref{eq:Phi} and components \eqref{eq:components}.
The following statements are equivalent:
\begin{enumerate}
\item $D$ is compact;
\item for some fixed $R>0$,
\begin{equation}\label{eq:local-condition}
 \Phi(a,a)\longrightarrow0,\qquad
 \omega_R(a)\longrightarrow0\qquad(|a|\to1);
\end{equation}
\item \eqref{eq:local-condition} holds for every fixed $R>0$;
\item all the following boundary limits hold as $|z|\to1$, with unrestricted approach:
\begin{equation}\label{eq:d-vanish}
 \Phi(z,z)=\sum_{j=1}^{N}\varphi_j(z)\psi_j(z)
 \longrightarrow0;
\end{equation}
\begin{equation}\label{eq:mixed-vanish}
 (1-|z|^2)^2\sum_{j=1}^{N}F_j'(z)\overline{K_j'(z)}
 \longrightarrow0;
\end{equation}
and, for every $i=1,\ldots,N$,
\begin{align}
 (1-|z|^2)^{i+1}\sum_{j=1}^{N}F_j^{(i)}(z)H_j'(z)
 &\longrightarrow0,\label{eq:hol-vanish}\\
 (1-|z|^2)^{i+1}\sum_{j=1}^{N}G_j^{(i)}(z)K_j'(z)
 &\longrightarrow0.\label{eq:anti-vanish}
\end{align}
\end{enumerate}
\end{theorem}

If all $\varphi_j,\psi_j$ extend continuously to $\overline{\D}$,
compactness is already characterized by
\cite[Theorem~3.2]{GZ}:
\[
 D\text{ is compact}
 \quad\Longleftrightarrow\quad
 \left.\sum_{j=1}^{N}\varphi_j\psi_j\right|_{\partial\D}=0.
\]
Indeed, its quotient map sends $T_{\alpha}$ to $\alpha|_{\partial\D}$.
In this case the centered condition $\omega_R(a)\to0$ is
automatic by uniform continuity. Theorem~\ref{thm:compact} also
covers bounded harmonic symbols without assuming continuous boundary values.

For every even $N\ge2$, the conditions of order $N$ in
Theorem~\ref{thm:compact}\textup{\,(iv)} cannot be omitted; see
Remark~\ref{rem:compact-even-order}.
In particular, no fixed derivative cutoff independent of the number
of product terms suffices in that criterion.

\section{Finite rank}\label{sec:finite}

In this section, we prove Theorem~\ref{thm:general} for the Toeplitz sum $D$ in
\eqref{eq:S}, using the components in \eqref{eq:centered}. We first introduce the following notation.

For a kernel $H(z,w)$ holomorphic in $z$ and antiholomorphic in $w$, we say that it has \emph{finite separation rank} if
\[
 H(z,w)=\sum_{j=1}^nA_j(z)\overline{B_j(w)},
 \quad A_j,B_j\in \mathrm{Hol}(\D).
\]
For a jointly holomorphic kernel, finite separation rank means a representation $H(z,w)=\sum_{j=1}^N A_j(z)B_j(w)$ with holomorphic factors. We use the convention appropriate to the holomorphy in the second variable.

\begin{remark}\label{rem:one-mixed}
Under conditions \textup{(i)} and \textup{(ii)} of
Theorem~\ref{thm:general}, the second mixed identity
\begin{equation}\label{eq:automatic-mixed}
\sum_{j=1}^{N}p_j(z)g_j(w)=0\qquad(z,w\in\D)
\end{equation}
is automatic. Applying $\partial\bar\partial$ to
\eqref{eq:full-symbol-sum} gives
\[
\sum_{j=1}^{N}\left[f_j'(z)q_j'(\bar z)
                  +p_j'(z)g_j'(\bar z)\right]=0.
\]
The first sum vanishes by \textup{(ii)}. Diagonal uniqueness
\cite[Lemma~3.1]{CKL} gives $\sum_jp_j'(z)g_j'(w)=0$.
Integrating from $0$ in both variables and using
$p_j(0)=g_j(0)=0$ gives \eqref{eq:automatic-mixed}.
The same argument works with the two mixed identities interchanged.
Equivalently, condition \textup{(ii)} can be checked on the full
symbols by
\begin{equation}\label{eq:full-symbol-derivative}
\sum_{j=1}^{N}(\partial_z\varphi_j)(z)(\partial_{\bar z}\psi_j)(z)=0
\qquad(z\in\D).
\end{equation}

\end{remark}

For later use, we recall a consequence of
\cite[Lemma~2.4 and proof of Proposition~2.5]{CHLZ}.
For $G=\sum_{\ell=1}^L \alpha_\ell\overline{\gamma_\ell}$ with
$\alpha_\ell,\gamma_\ell\in H^2(\D)$, set
\[
\mathcal W_G(z,w)=\int_\D
\frac{G(\xi)}{(1-z\bar\xi)^2(1-\bar w\xi)^2}\,dA(\xi).
\]
Then
\begin{equation}\label{eq:weighted-rigidity}
(1-z\bar w)^4\mathcal W_G(z,w)
\text{ has finite separation rank}
\quad\Longrightarrow\quad \partial_z\partial_{\bar z} G=0.
\end{equation}
Indeed, the cited weighted identity and finite-rank argument apply to
this finite sum and give a finite-rank analytic Toeplitz operator
with bounded symbol
$\rho(\xi)=|\xi|^2(1-|\xi|^2)^4
\partial_\xi\partial_{\bar\xi}G(\xi)$;
Luecking's theorem \cite{Luecking} gives $\rho=0$, and hence
$\partial_z\partial_{\bar z}G=0$.

The following finite-sum version of the argument in
\cite[Proposition~2.5]{CHLZ} will be needed; the cross terms
in its analytic compression are retained explicitly.

\begin{proposition}\label{prop:input}
For $D$ and its  components in \eqref{eq:S} and
\eqref{eq:centered}, finite rank of $D$ implies
\[
\sum_{j=1}^{N}\varphi_j\psi_j=0,\qquad
\sum_{j=1}^{N}p_j'(z)g_j'(\bar z)=0,\qquad
\sum_{j=1}^{N}f_j'(z)q_j'(\bar z)=0.
\]
Moreover,
\[
\sum_{j=1}^{N}p_j(z)g_j(w)=0,\qquad
\sum_{j=1}^{N}f_j(z)q_j(w)=0\qquad(z,w\in\D).
\]
\end{proposition}

\begin{proof}
The analytic components belong to $H^2(\D)$.
The analytic compression of $D$ has the kernel
$\mathcal F(z,w)=(PDK_w)(z)$, which has finite separation rank.
Block multiplication relative to $\bh=\A\oplus\overline{\Azer}$ gives
\begin{align}
\mathcal F(z,w)
={}&\mathcal W_{\sum_jp_jg_j(\bar\cdot)}(z,w)
+\sum_{j=1}^N\left\langle
\bar P_0\bigl(q_j(\bar\cdot)K_w\bigr),
\bar P_0(\overline{f_j}K_z)\right\rangle\notag\\
&+K_w(z)\sum_{j=1}^N\Bigl[
(\varphi_j(0)+f_j(z))
(\psi_j(0)+p_j(z)+q_j(\bar w))\\
&\hspace{38mm}+g_j(\bar w)(\psi_j(0)+q_j(\bar w))\Bigr]\notag.
\end{align}
For general bounded harmonic symbols, obtain this identity by
common radial dilation of the full symbols. The Toeplitz factors are
uniformly bounded and converge strongly, their analytic components
converge in $H^2(\D)$, the products defining $\mathcal W$ converge in
$L^1(\D, dA)$, and the projected products with fixed kernels converge in $L^2(\D, dA)$.
Thus the bounded-component calculation passes to the limit.

Multiplication by $(1-z\bar w)^4$ makes the last two lines finite
separated sums. The projected inner products have the same property
by \cite[Lemma~2.3]{CHLZ}. Hence
$(1-z\bar w)^4\mathcal W_{\sum_{j=1}^Np_jg_j(\bar\cdot)}$ has finite
separation rank. Applying \eqref{eq:weighted-rigidity}, with
$g_j(\bar\xi)=\overline{\widetilde g_j(\xi)}$ and
$\widetilde g_j(\xi)=\overline{g_j(\bar\xi)}\in H^2(\D)$, gives
$\sum_{j=1}^Np_j'(z)g_j'(\bar z)=0$.
Since the harmonic Bergman projection commutes with pointwise complex
conjugation, for every $u\in\bh$ we have
\[
 \left(\sum_{j=1}^N
 T_{\overline{\varphi_j}}T_{\overline{\psi_j}}\right)u
 =\overline{D(\overline{u})}.
\]
Thus the range of this operator is
$\{\overline{v}:v\in\operatorname{Ran}D\}$, so it has the same finite
rank as $D$. Applying the preceding argument to the conjugated symbols
and taking complex conjugates yields
$\sum_{j=1}^N f_j'(z)q_j'(\bar z)=0$.
Diagonal uniqueness \cite[Lemma~3.1]{CKL}, followed by integration
from $0$ in both variables and the centering of the components,
gives both two-variable identities.

These identities cancel the mixed terms in the analytic compression.
For the projected cross terms, fix $z,w$ and write
\[
 B_{z,w}(f,q)=\big\langle\bar P_0((Uq)K_w),
                            \bar P_0(\overline f K_z)\big\rangle .
\]
The form $B_{z,w}$  is complex bilinear in $f,q$. Expanding the $f_j$ in a basis of their span,
the identity $\sum_j f_j\otimes q_j=0$ gives $\sum_j B_{z,w}(f_j,q_j)=0$.
Both two-variable identities persist under common radial dilation, so for $0<r<1$,
\[
P\Bigl(\sum_{j=1}^NT_{\varphi_j(r\cdot)}T_{\psi_j(r\cdot)}\Bigr)\big|_{\A}
=\TT_{\sum_{j=1}^N\varphi_j(r\cdot)\psi_j(r\cdot)}.
\]
Both sides converge strongly, and therefore
$PD|_{\A}=\TT_{\sum_{j=1}^N\varphi_j\psi_j}$.
Its symbol is bounded. Luecking's theorem \cite{Luecking} applied to
this finite-rank analytic Toeplitz operator yields
$\sum_{j=1}^N\varphi_j\psi_j=0$ almost everywhere, and hence everywhere by
continuity.
\end{proof}

For a bounded operator $S:\overline{\Azer}\to\A$, we let
\[
\ell_w(\zeta)=\frac{1}{(1-\bar\zeta w)^2}-1,
\qquad
\mathscr K_S(z,w)=(S\ell_w)(z).
\]
Thus $\ell_w$ is the reproducing kernel of $\overline{\Azer}$,
and $\mathscr K_S$ is holomorphic in both $z$ and $w$.
The linear span of the functions $\ell_w$, $w\in\D$, is dense in
$\overline{\Azer}$.

For arbitrary $f,h\in\operatorname{Hol}(\D)$, define the holomorphic kernel
\[
\mathscr K_{f,h}(z,w)
:=\frac{2zw}{(z-w)^3}\mathscr R_{f,h}(z,w).
\]
Here and below the quotient is extended across the diagonal:
by \eqref{eq:Rintegral}, $\frac{\mathscr R_{f,h}(z,w)}{(z-w)^3}$ is
holomorphic on $\D^2$ and has diagonal value $-f'(z)h'(z)/6$.
This definition does not presume boundedness of $f,h$ or of an
individual semicommutator associated with them.

When $f,h\in H^\infty(\D)$, put
$(T_f,T_h]=T_fT_h-T_{fh}$.
The projection formula \cite[Eq.~(2.30)]{CHLZ}, together with
the coanalytic projection identity in the proof of
\cite[Theorem~2.8]{CHLZ} gives, for $\alpha\in\{f,h,fh\}$,
\begin{align}
T_{\alpha}\ell_w&=J_{\alpha}(\cdot,w)+\alpha(w)\ell_w,\\
J_{\alpha}(z,w)&=
\frac{z^2\alpha(z)+w^2\alpha(w)}{(z-w)^2}
-\frac{2zw(F_{\alpha}(z)-F_{\alpha}(w))}{(z-w)^3}-\alpha(z).
\end{align}

Consequently, $(T_f,T_h]$ vanishes on $\A$, and its cross kernel is
\begin{equation}\label{eq:semicommutator-cross-kernel}
\begin{aligned}
\mathscr K_{f,h}(z,w)
&=((T_f,T_h]\ell_w)(z)\\
&=f(z)J_h(z,w)+h(w)J_f(z,w)-J_{fh}(z,w)\\
&=\frac{2zw}{(z-w)^3}\mathscr R_{f,h}(z,w).
\end{aligned}
\end{equation}
The factor $z$ and density of the kernel span show that
$(T_f,T_h]$ maps $\overline{\Azer}$ into $\Azer$.

\begin{lemma}
\label{lem:cancellation}
For the symbols and  decompositions in
\eqref{eq:S} and \eqref{eq:centered}, suppose
\[
\sum_{j=1}^N \varphi_j(z)\psi_j(z)=0\qquad(z\in\D)
\]
and
\[
\sum_{j=1}^N f_j(z)q_j(w)=0,
\qquad
\sum_{j=1}^N p_j(z)g_j(w)=0
\qquad(z,w\in\D).
\]
There are bounded operators
$B_+,B_-:\overline{\Azer}\to\Azer$, extended by zero on
$\A$, such that
\[
D=B_++UB_-U.
\]
Their cross kernels are
\[
\mathscr K_+(z,w)=\sum_{j=1}^N\mathscr K_{f_j,p_j}(z,w),
\qquad
\mathscr K_-(z,w)=\sum_{j=1}^N\mathscr K_{g_j,q_j}(z,w).
\]
If all the analytic components are bounded, then
\[
B_+=\sum_{j=1}^N (T_{f_j},T_{p_j}],\qquad
B_-=\sum_{j=1}^N (T_{g_j},T_{q_j}].
\]
In general, these sums denote the combined bounded blocks
obtained by common radial dilation; they do not assert
boundedness of the individual semicommutators.

The operator $D$ has finite rank if and only if both
$\mathscr K_+$ and $\mathscr K_-$ have finite separation rank.
In that case,
\[
\rank D=\rank B_++\rank B_-.
\]
For each sign, the rank of $B_\pm$ equals the least number of
terms in a separated representation of $\mathscr K_\pm$.
\end{lemma}

\begin{proof}
First suppose that all analytic components are bounded.
The constant-basis argument in \cite[proof of Theorem~3.2]{CKL}
has the following bilinear form: if $\sum_{j=1}^Np_j(z)q_j(w)=0$,
then $\sum_{j=1}^N\mathcal B(p_j,q_j)=0$ for every bilinear map
$\mathcal B$ on the two finite spans. Indeed, writing
$p_j=\sum_\nu a_{\nu j}p^{(\nu)}$ in a basis gives
$\sum_{j=1}^Na_{\nu j}q_j=0$ for every $\nu$.

Apply this observation to the
bilinear maps $(f,q)\mapsto T_fT_{Uq}$ and
$(p,g)\mapsto T_{Ug}T_p$. The hypotheses give
\[
\sum_{j=1}^NT_{f_j}T_{Uq_j}=0,\qquad
\sum_{j=1}^NT_{Ug_j}T_{p_j}=0.
\]
They also give
$\sum_{j=1}^Nf_j(Uq_j)=\sum_{j=1}^N(Ug_j)p_j=0$ by taking $w=\bar z$.
Expanding the products, with all constant terms retained,
therefore yields
\[
\begin{split}
D-\sum_{j=1}^NT_{\varphi_j\psi_j}
={}&\sum_{j=1}^N (T_{f_j},T_{p_j}]
 +U\left(\sum_{j=1}^N (T_{g_j},T_{q_j}]\right)U\\
&+\sum_{j=1}^N\bigl(T_{f_j}T_{Uq_j}-T_{f_j(Uq_j)}\bigr)\\
&+\sum_{j=1}^N\bigl(T_{Ug_j}T_{p_j}-T_{(Ug_j)p_j}\bigr).
\end{split}
\]
The last two sums are zero, and $\sum_{j=1}^NT_{\varphi_j\psi_j}=0$.
The asserted identity and kernel formulas now follow from
\eqref{eq:semicommutator-cross-kernel} and $UT_qU=T_{Uq}$.

For arbitrary bounded harmonic symbols, use
$\varphi_{j,r}(z)=\varphi_j(rz)$ and $\psi_{j,r}(z)=\psi_j(rz)$,
with a common dilation parameter $0<r<1$ for all symbols.
All three cancellation identities persist, so the bounded-component
case applied to $D_r=\sum_jT_{\varphi_{j,r}}T_{\psi_{j,r}}$
gives $D_r=B_{+,r}+UB_{-,r}U$, where
\[
B_{+,r}=P_0D_r\bar P_0,\qquad
B_{-,r}=U\bar P_0D_rP_0U.
\]
As in \cite[proof of Proposition~2.5]{CHLZ}, dominated convergence
for the bounded full symbols gives strong convergence of their
Toeplitz operators. Their norms are uniformly bounded, so
$D_r\to D$ strongly and
$\|D_r\|\le\sum_j\|\varphi_j\|_\infty\|\psi_j\|_\infty$.
The fixed compressions therefore converge strongly to
$B_+=P_0D\bar P_0$ and $B_-=U\bar P_0DP_0U$, proving the asserted
bounded decomposition.
The analytic components and their primitives converge locally
uniformly. Evaluating $B_{\pm,r}\ell_w$ and using the continuity of
point evaluation on $\Azer$ identifies the limiting kernels with
$\mathscr K_\pm$ through \eqref{eq:semicommutator-cross-kernel}.
Thus only the combined blocks are asserted to be bounded; no bounded
operator is assigned to an individual unbounded analytic component.

Finally, a bounded cross operator $S$ has finite rank if and
only if $\mathscr K_S$ has finite separation rank. Indeed,
a finite-dimensional range gives a separated representation
by choosing a basis of that range. Conversely, a finite
separated representation places all vectors $S\ell_w$ in a
fixed finite-dimensional space. Their linear span is then
closed, and density of the span of $\ell_w$ implies that it
contains the entire range of $S$. This also proves equality
between the operator rank and the minimum number of
separated terms. In the decomposition of $D$, the two
blocks have orthogonal ranges and act on different input
summands. Their ranks therefore add, proving the final
assertions.
\end{proof}

All spans over rational function fields below are taken inside
the corresponding fields of meromorphic functions.

\begin{lemma}\label{lem:quotient}
Let $N_0\ge1$ be an integer, and let
\[
\mathcal{R}(z,w)=\sum_{\ell=1}^{N}\alpha_\ell(z)\gamma_\ell(w),
\qquad \alpha_\ell,\gamma_\ell\in\operatorname{Hol}(\D).
\]
Assume that the function
\[
K(z,w)=\frac{\mathcal{R}(z,w)}{(z-w)^{N_0}},\qquad z\ne w,
\]
extends holomorphically to $\D^2$. Choose bases $\xi_1,\ldots,\xi_r$ and
$\eta_1,\ldots,\eta_m$ for
\[
\operatorname{span}_{\C(z)}\{\alpha_1,\ldots,\alpha_N\},
\qquad
\operatorname{span}_{\C(w)}\{\gamma_1,\ldots,\gamma_N\},
\]
respectively. These bases may be chosen from the original families,
and we make this choice so that all basis functions are holomorphic.
After clearing the rational coefficients, write
\begin{equation}\label{eq:quotient-representation}
\mathcal{R}(z,w)=\frac{1}{A(z)B(w)}
\sum_{i=1}^{r}\sum_{j=1}^{m}
P_{ij}(z,w)\xi_i(z)\eta_j(w),
\end{equation}
where $A,B$ are nonzero one-variable polynomials and
$P_{ij}\in\C[z,w]$. Then $K$ has finite separation rank if and only if
\begin{equation}\label{eq:divisibility}
(z-w)^{N_0}\mid P_{ij}(z,w)
\qquad(1\le i\le r,\ 1\le j\le m).
\end{equation}
When \eqref{eq:divisibility} holds, put
$Q_{ij}=\frac{P_{ij}}{(z-w)^{N_0}}$ and expand
\begin{equation}\label{eq:quotient-coefficients}
Q_{ij}(z,w)=
\sum_{s,t\ge0}
c_{is,jt}z^s w^t.
\end{equation}
Let $M=(c_{is,jt})_{(i,s),(j,t)}$ be the
finite coefficient matrix, using common finite ranges for $s$
and $t$ and inserting zero coefficients where necessary. Then
\begin{equation}\label{eq:quotient-rank}
\operatorname{srank}K=\rank M.
\end{equation}
Both the divisibility criterion and the resulting rank are independent
of the finite separated representation of $\mathcal{R}$, the chosen bases, and
the common denominators.
\end{lemma}

\begin{proof}
We first establish that the products $\xi_i(z)\eta_j(w)$ are linearly
independent over $\C(z,w)$. A rational-coefficient relation can be
multiplied by a common denominator to give a polynomial-coefficient
relation. Expanding its coefficients gives a finite relation among
\[
z^s \xi_i(z)\,w^t \eta_j(w).
\]
Every finite family of distinct functions $z^s \xi_i(z)$ is
linearly independent over $\C$: a relation among them would be a
polynomial-coefficient relation among the $\xi_i$. The same statement
holds for $w^t \eta_j(w)$. First fixing $w$ and using independence
in $z$, and then using independence in $w$, shows that every
coefficient in the expanded relation is zero. This proves the
claimed independence over $\C(z,w)$.
It also treats the zero-kernel case: if $K=0$, all the $P_{ij}$ are
zero and every assertion follows. We henceforth assume $K\ne0$.

Suppose that $\operatorname{srank}K=n<\infty$. In a minimal representation,
the left and right factors are separately linearly independent. Their
evaluation vectors span the respective coordinate spaces, since a linear
functional annihilating every evaluation vector gives a relation among
the factors. Choose $w_1,\ldots,w_n$ so that the slices
$\delta_\ell=K(\cdot,w_\ell)$, $1\le \ell\le n$, form a basis of the
column space. Choose $z_1,\ldots,z_n$ so that the matrix
$V=(K(z_i,w_j))_{i,j=1}^n$ is invertible. Then
\[
 K(z,w)=\sum_{\ell=1}^n \delta_\ell(z)\epsilon_\ell(w),
 \qquad
 \begin{pmatrix} \epsilon_1(w)\\ \vdots\\ \epsilon_n(w) \end{pmatrix}
 =
 V^{-1}
 \begin{pmatrix} K(z_1,w)\\ \vdots\\ K(z_n,w) \end{pmatrix}.
\]
Since $K=\mathcal{R}/(z-w)^{N_0}$ and $\mathcal{R}=\sum_k \alpha_k(z)\gamma_k(w)$,
evaluating at $w=w_\ell$ and dividing by $(z-w_\ell)^{N_0}$
shows that $\delta_\ell$ lies in the $\C(z)$-span of the $\xi_i$.
Likewise each row slice, and hence each $\epsilon_\ell$, lies in the
$\C(w)$-span of the $\eta_j$.
Expanding the two families in these bases, the coefficient of
$\xi_i(z)\eta_j(w)$ in $K$ is a finite sum of products of one-variable
rational functions. Clearing their denominators gives a common
denominator $A_1(z)B_1(w)$. By the independence proved above and
\eqref{eq:quotient-representation},
\[
 \frac{P_{ij}(z,w)}{A(z)B(w)(z-w)^{N_0}}
 =\frac{N_{ij}(z,w)}{A_1(z)B_1(w)},\qquad N_{ij}\in\C[z,w].
\]
As $z-w$ is prime and divides neither $A_1(z)$ nor $B_1(w)$,
cross multiplication gives \eqref{eq:divisibility}.

Conversely, suppose that \eqref{eq:divisibility} holds. Equations
\eqref{eq:quotient-representation} and
\eqref{eq:quotient-coefficients} give
\begin{equation}\label{eq:quotient-meromorphic-factors}
K(z,w)=
\sum_{i,s,j,t}
c_{is,jt}
\frac{z^s \xi_i(z)}{A(z)}
\frac{w^t \eta_j(w)}{B(w)}.
\end{equation}
This is a finite separated representation by meromorphic factors.
Choose nonempty open sets $U,V\subset\D$ avoiding the zeros of
$A$ and $B$, respectively. On these sets the factors in
\eqref{eq:quotient-meromorphic-factors} are holomorphic, and each
of the two finite families of factors is linearly independent.
Indeed, a relation on an open set extends as a meromorphic identity,
and multiplication by $A$ or $B$ reduces it to the independence
used at the beginning of the proof.

It follows that the separation rank on $U\times V$ is
$n=\rank M$. To see the lower bound explicitly, the evaluation
vectors of the factors in $w$ span their full coordinate space.
Their images under $M$ therefore span its range. Independence
of the factors in $z$ shows that the column-slice space has
dimension $\rank M$. The reverse inequality follows by
factoring the finite matrix $M$.

To recover holomorphic factors on all of $\D$, choose
$w_1,\ldots,w_n\in V$ and $z_1,\ldots,z_n\in U$ giving an
invertible evaluation matrix for this rank-$n$ kernel. As above,
set $\delta_\ell(z)=K(z,w_\ell)$ and obtain the column vector
$(\epsilon_\ell(w))_{\ell=1}^n$ by applying the inverse evaluation matrix
to the column vector of row slices $(K(z_i,w))_{i=1}^n$.
These functions are holomorphic on all of $\D$
because $K$ is holomorphic on $\D^2$. The resulting identity
\[
K(z,w)=\sum_{\ell=1}^{n}\delta_\ell(z)\epsilon_\ell(w)
\]
holds on $U\times V$. Applying the one-variable identity theorem
successively in $z$ and $w$ extends it to $\D^2$.
Thus the global separation rank is at most $n$, while
restriction to $U\times V$ gives the opposite inequality.
This proves both sufficiency and \eqref{eq:quotient-rank}.
Finally, any other finite separated representation of $\mathcal{R}$, choice
of rational bases, or choice of common denominators gives the same
test and the same rank, since both have been identified with
intrinsic properties of $K$.
\end{proof}

\begin{proof}[Proof of Theorem~\ref{thm:general}]
If $D$ has finite rank, Proposition~\ref{prop:input} gives
\eqref{eq:full-symbol-sum} and \eqref{eq:mixed-identities}.
Under these two identities,
Remark~\ref{rem:one-mixed} gives \eqref{eq:automatic-mixed}, and
Lemma~\ref{lem:cancellation} yields
\[
D=B_++UB_-U,\qquad
\mathscr K_\pm(z,w)=\frac{2zw}{(z-w)^3}\mathcal{R}_\pm(z,w).
\]
For bounded analytic components, these blocks are
\begin{equation}\label{eq:direct-blocks}
B_+=\sum_{j=1}^{N} (T_{f_j},T_{p_j}],\qquad
B_-=\sum_{j=1}^{N} (T_{g_j},T_{q_j}].
\end{equation}
For general bounded full symbols, they are the combined bounded
blocks obtained in Lemma~\ref{lem:cancellation} by common radial
dilation.

Each block has the separation rank of its cross kernel.
Multiplying by the nonzero one-variable factors $2z$ and $w$
preserves separation rank, first away from their zeros and then on
the whole disk by the holomorphic slice construction in
Lemma~\ref{lem:quotient}. Applying that lemma with $N_0=3$ gives
\[
B_\pm\text{ has finite rank}
\quad\Longleftrightarrow\quad
\mathcal{R}_\pm\text{ satisfies }\eqref{eq:rd3},
\qquad \rank B_\pm=\rank M_{\mathcal{R}_\pm}
\]
whenever these conditions hold. This proves the characterization.
The two cross blocks have orthogonal ranges and act on different
input summands, so their ranks add, giving \eqref{eq:rank-general}.
They annihilate constants and exchange the zero-constant analytic
and coanalytic summands. With $C_\pm=B_\pm U|_{\Azer}$, this proves
\eqref{eq:general-diagonal-blocks}; equivalently,
$DU\simeq0_\C\oplus C_+\oplus C_-$.

The family in \eqref{eq:two-orders-rankr-def} realizes every positive
rank by \eqref{eq:two-orders-rankr-range}, and the zero operator realizes rank zero.
\end{proof}

\begin{remark}\label{rem:general-primitives}
Under \eqref{eq:full-symbol-sum} and \eqref{eq:mixed-identities},
all product primitives in the combined numerators can be eliminated.
The two mixed identities cancel the mixed terms in the full symbol
sum. Its centered analytic and coanalytic parts therefore give
\[
\begin{aligned}
\sum_{j=1}^{N}\bigl[f_jp_j+\varphi_j(0)p_j+\psi_j(0)f_j\bigr]&=0,\\
\sum_{j=1}^{N}\bigl[g_jq_j+\varphi_j(0)q_j+\psi_j(0)g_j\bigr]&=0.
\end{aligned}
\]
Taking primitives that vanish at zero yields
\[
\begin{aligned}
\sum_{j=1}^{N}F_{f_jp_j}
&=-\sum_{j=1}^{N}\bigl[\varphi_j(0)F_{p_j}+\psi_j(0)F_{f_j}\bigr],\\
\sum_{j=1}^{N}F_{g_jq_j}
&=-\sum_{j=1}^{N}\bigl[\varphi_j(0)F_{q_j}+\psi_j(0)F_{g_j}\bigr].
\end{aligned}
\]
Thus the numerators $\mathcal{R}_+$ and $\mathcal{R}_-$ can be expressed using only
primitives of the individual analytic components. The polynomial
divisibility condition \eqref{eq:rd3} still applies to those combined
expressions.
\end{remark}

The following monomial example of Dong and Zhou
\cite[Theorem~5.4 and Corollary~5.6(a), arXiv version]{DZ} illustrates the rank formula
in Theorem~\ref{thm:general}. For $r\ge1$, set
\begin{equation}\label{eq:two-orders-rankr-def}
 D_r=T_zT_{z^r}-T_{z^{r+1}}T_1 .
\end{equation}
With $k_1=1$ and $k_2=r$, their rank formula gives
$\max\{0,r-1,r\}=r$. Directly, the monomial projection formula
\cite[Lemma~3.1]{DZ} gives $D_r|_{\A}=0$ and
$D_r\bar z^k=-kz^{r+1-k}/[(r+1)(r+2)]$ for $1\le k\le r$,
while $D_r\bar z^k=0$ for $k\ge r+1$. Thus
\begin{equation}\label{eq:two-orders-rankr-range}
 \operatorname{Ran}D_r=\operatorname{span}\{z,\ldots,z^r\},\qquad
 \rank D_r=r.
\end{equation}
In particular, odd ranks occur for two-term sums, whereas a Toeplitz
commutator on $b^2(\D)$ cannot have odd rank \cite{ChenKooLee}.

To see how \eqref{eq:rd3} distinguishes compactness from finite rank,
consider $D=T_zT_{e^z}-T_{ze^z}T_1$. Its diagonal symbol sum is zero.
Here $\mathcal{R}_-=0$ and $\mathcal{R}_+=\mathscr R_{z,e^z}$, since constants drop out
of \eqref{eq:Rintegral}. Writing $t=z-w$ in \eqref{eq:Rdef} gives
\[
 \mathcal{R}_+(z,w)=e^w(1+t+t^2/2)-e^z.
\]
The quotient $\mathcal{R}_+/t^3$ extends holomorphically across the diagonal.
In the rational bases $\{1,e^z\}$ and $\{1,e^w\}$, however, the
coefficient of $e^z\cdot1$ is $-1$, so \eqref{eq:rd3} fails and
$D$ has infinite rank. This conclusion also follows indirectly from
\cite[Theorem~2.8]{CHLZ}: transposition of the cross kernel interchanges
$(T_z,T_{e^z}]$ and $(T_{e^z},T_z]$, so finite rank of the former
would imply finite rank of $[T_z,T_{e^z}]$, contrary to the cited
classification since $1,z,e^z$ are linearly independent and $e^z$
has no rational primitive.
All symbols extend continuously to the closed disk, so compactness
already follows from the continuous-symbol Toeplitz algebra theorem
\cite[Theorem~3.2]{GZ}. It is also immediate from the weighted limits
in Theorem~\ref{thm:compact}\textup{\,(iv)}, because all derivatives
are bounded there. Thus the example illustrates how the cubic divisibility
test detects infinite rank among compact continuous-symbol semicommutators.

\section{Compactness}\label{sec:compact}

This section proves Theorem~\ref{thm:compact} for the Toeplitz sum $D$
in \eqref{eq:S}. We first characterize vanishing finite sums of
Toeplitz products on the half-plane and then pass to boundary limits.

The criterion retains all cancellations between the summands.
The centered expression in \eqref{eq:omega} is
\[
 \sum_{j=1}^{N}[\varphi_j(z)-\varphi_j(a)]
                  [\psi_j(w)-\psi_j(a)].
\]
Theorem~\ref{thm:compact}\textup{\,(iv)} gives $2N+2$ scalar limits
and differentiates each symbol component at most $N$ times. The second mixed condition
\[
 (1-|z|^2)^2\sum_{j=1}^{N}H_j'(z)\overline{G_j'(z)}
 \longrightarrow0
\]
follows from the listed conditions, as the proof will show.

We use the localization theorem of Liu and Lin
\cite[Theorem~1.1(2) and Corollary~1.3]{LiuLin}. For a finite
sum $T$ of products of bounded-symbol Toeplitz operators it gives
\begin{equation}\label{eq:localization}
 T\text{ is compact}\quad\Longleftrightarrow\quad
 \exists R>0:\quad
 \lim_{|a|\to1}\sup_{b\in D_\beta(a,R)}
 |\langle Tr_b,r_a\rangle|=0,
\end{equation}
where $r_a$ is the normalized harmonic reproducing kernel.
Their theorem is stated for real harmonic spaces. We use its
complexification, as explained in \cite[Section~3, preceding
Lemma~3.4]{CHLZ}, by applying the real result to the real and
imaginary parts of the operator polynomial.

\subsection{Toeplitz products on the half-plane}\label{sec:halfplane}

Let $\HH=\{z\in\C:\operatorname{Re}z>0\}$. Let
$Q_{\HH}:L^2(\HH, dA)\to b^2(\HH)$ be the orthogonal projection.
For  $u \in L^\infty(\HH, dA) $, the Toeplitz operator is defined by
\[T_u^{\HH}=Q_{\HH}M_u|_{b^2(\HH)}.\]
We establish the following characterization.

\begin{proposition}
\label{prop:halfplane-zero}
Let $u_j,v_j\in h^\infty(\HH)$, $1\le j\le N$. Then
\begin{equation}\label{eq:halfplane-zero}
 \sum_{j=1}^{N}T_{u_j}^{\HH}T_{v_j}^{\HH}=0
\end{equation}
if and only if there exists $\phi\in h^\infty(\HH)$ such that
\begin{equation}\label{eq:additive}
 \sum_{j=1}^{N}u_j(z)v_j(w)=\phi(z)-\phi(w),\qquad z,w\in\HH.
\end{equation}
\end{proposition}

We apply \cite[Lemma~3.2 and its proof]{CHLZ} in the following
$L^2(\D, dA)$ form. For $\delta,\epsilon\in L^2(\D, dA)$, put
\[
\mathcal C_{\delta,\epsilon}(z,w)=
\left\langle\bar P(\bar\epsilon K_w),\bar P(\bar\delta K_z)\right\rangle.
\]
Then
\begin{equation}\label{eq:weighted-cross-kernel}
(1-z\bar w)^4\mathcal C_{\delta,\epsilon}(z,w)
\quad\text{has finite separation rank},
\end{equation}
with a common finite bound independent of $\delta,\epsilon$.
For bounded $\delta,\epsilon$, the full coanalytic projection adds only the
constant component to the kernel in the cited lemma:
\[
\mathcal C_{\delta,\epsilon}(z,w)
=(N_\delta N_\epsilon^*K_w)(z)+(P\delta)(z)\overline{(P\epsilon)(w)},
\qquad N_\delta=PM_\delta|_{\overline{zA^2(\D)}}.
\]
For clarity, one uniform rank bound is $n_0=35\cdot5^4$.
For the full projection the constant term cancels the subtraction in
\cite[Eq.~(3.3)]{CHLZ}, leaving a double integral. The fourth power
in \cite[Eq.~(3.4)]{CHLZ} has at most $\binom{7}{3}=35$ terms.
In each term a squared kernel denominator cancels. If it is an outer
denominator, expansion in $z$ or $\bar w$ gives at most five separated
terms. If it is the middle denominator, the remaining denominators
separate and the numerator has degree at most four in each of
$z,\bar w,\bar\xi,\eta$, hence at most $5^4$ monomials.
For $\delta,\epsilon\in L^2(\D,dA)$, approximate both in $L^2(\D,dA)$ by bounded functions.
Since $K_z,K_w$ are bounded for fixed $z,w$, the kernels converge
pointwise. Every $(n_0+1)\times(n_0+1)$ evaluation determinant is zero before
and after passage to the limit, proving \eqref{eq:weighted-cross-kernel}.

There is an orthogonal decomposition
\begin{equation}\label{eq:H-decomposition}
 b^2(\HH)=A^2(\HH)\oplus\overline{A^2(\HH)}.
\end{equation}
With $dA=\pi^{-1}dx\,dy$, the analytic reproducing kernel of $A^2(\HH)$ is
\begin{equation}\label{eq:H-kernel}
 K_w^{\HH}(z)=\frac{1}{(z+\bar w)^2},\quad z, w \in \HH.
\end{equation}

Let $P_{\HH}:L^2(\HH)\to A^2(\HH)$ be the analytic projection.
For $v\in L^\infty(\HH, dA)$, define the Toeplitz operator
\[\TT_v^{\HH}=P_{\HH}M_v|_{A^2(\HH)}.\]
Let $\bar P_{\HH}$ denote the orthogonal projection of $L^2(\HH)$
onto $\overline{A^2(\HH)}$. For $\alpha\in H^\infty(\HH)$, the block
matrix of $T_\alpha^{\HH}$ relative to \eqref{eq:H-decomposition} is
\begin{equation}\label{eq:block}
\begin{aligned}
 T_\alpha^{\HH}
 &=
 \begin{pmatrix}
   M_\alpha|_{A^2(\HH)} & N_\alpha \\
   0 & L_\alpha
 \end{pmatrix},\\
 N_\alpha
 &=P_{\HH}M_\alpha|_{\overline{A^2(\HH)}},
 \qquad
 L_\alpha=\bar P_{\HH}M_\alpha|_{\overline{A^2(\HH)}}.
\end{aligned}
\end{equation}
Since multiplication by $\bar\alpha$ preserves $\overline{A^2(\HH)}$,
we have $L_\alpha^*=M_{\bar\alpha}|_{\overline{A^2(\HH)}}$.
Thus, for $\alpha,\gamma\in H^\infty(\HH)$,
\[
 (L_\alpha L_\gamma)^*
 =L_\gamma^*L_\alpha^*
 =M_{\overline{\alpha\gamma}}|_{\overline{A^2(\HH)}}
 =L_{\alpha\gamma}^*,
\]
and hence $L_\alpha L_\gamma=L_{\alpha\gamma}$.

The next formula is the half-plane analogue of the disk projection identity in \cite[Eq.~(2.30)]{CHLZ}.

\begin{lemma}
\label{lem:H-cross}
Let $\alpha\in H^\infty(\HH)$ and let $A'=\alpha$. Then
\begin{equation}\label{eq:C-a}
 (N_\alpha\overline{K_w^{\HH}})(z)
 =C_\alpha(z,w)=
 \frac{(z-w)[\alpha(z)+\alpha(w)]-2[A(z)-A(w)]}{(z-w)^3}.
\end{equation}
The right side has a removable singularity on the diagonal.
\end{lemma}
\begin{proof}
For $\varepsilon>0$, put
$\alpha_\varepsilon(\xi)= \frac{\alpha(\xi)}{(1+\varepsilon\xi)^2}$.
Since $\operatorname{Re}\xi>0$, we have
$|\alpha_\varepsilon|\le\|\alpha\|_\infty$, and
$\alpha_\varepsilon\in A^2(\HH)$. Its cross kernel is
\begin{equation}\label{eq:H-integral}
 \int_{\HH}
 \frac{\alpha_\varepsilon(\xi)}{(\bar\xi+z)^2(\bar\xi+w)^2}\,dA(\xi).
\end{equation}
For $z\ne w$,
\begin{align*}
 \frac1{(x+z)^2(x+w)^2}
 ={}&\frac1{(z-w)^2}
       \left(\frac1{(x+z)^2}+\frac1{(x+w)^2}\right)\\
 &+\frac2{(z-w)^3}
       \left(\frac1{x+z}-\frac1{x+w}\right).
\end{align*}
The first two integrals are evaluated by the analytic reproducing
formula. Define
\[
 I_\varepsilon(z)=\int_{\HH}
                  \frac{\alpha_\varepsilon(\xi)}{\bar\xi+z}\,dA(\xi).
\]
The integrand is $O(|\xi|^{-3})$ at infinity, and the denominator
is bounded away from zero on bounded regions. Thus the integral
is absolutely convergent and can be differentiated on compact
parameter sets. The reproducing formula gives
$I_\varepsilon'(z)=-\alpha_\varepsilon(z)$, and hence
\[
 I_\varepsilon(z)-I_\varepsilon(w)
 =-\int_w^z\alpha_\varepsilon(\xi)\,d\xi.
\]
Substitution into \eqref{eq:H-integral} proves \eqref{eq:C-a}
for $\alpha_\varepsilon$.

Let $\varepsilon\rightarrow0$. The integrand in
\eqref{eq:H-integral} is dominated by a constant times the product
of two square-integrable kernel functions, while the primitive
differences on the right converge locally uniformly. The formula
therefore holds for $\alpha$. Taylor expansion shows that its numerator
vanishes to order at least three at $z=w$, proving removability.
\end{proof}

\begin{proof}[Proof of Proposition~\ref{prop:halfplane-zero}]
We first prove necessity. Write
$u_j=f_j+\overline{g_j}$ and $v_j=h_j+\overline{k_j}$.
Let $\chi(\zeta)=\frac{1+\zeta}{1-\zeta}$ be the Cayley map.
As in \cite[proof of Lemma~3.3]{CHLZ}, use the unitary maps
on the analytic and coanalytic spaces, respectively,
\begin{equation}\label{eq:Cayley-unitary}
 U_+F=\chi'(F\circ\chi),\qquad
 U_-\overline F=\overline{\chi'}(\overline F\circ\chi).
\end{equation}
Put $F_j=f_j\circ\chi$, $G_j=g_j\circ\chi$,
$H_j=h_j\circ\chi$, and $K_j=k_j\circ\chi$.
After adjusting constants so that $G_j(0)=K_j(0)=0$, the
harmonic decomposition  gives
$F_j,G_j,H_j,K_j\in H^2(\D)$, since the pulled-back full
symbols are bounded and harmonic.

Let
\[
 \sigma(\zeta)=\frac{\chi'(\zeta)}{\overline{\chi'(\zeta)}},
 \qquad |\sigma|=1.
\]
First fix $0<r<1$. For each disk component $X$ write
$X_r(\zeta)=X(r\zeta)$ and pull the corresponding full symbols back to
the half-plane by $\chi^{-1}$. All four dilated analytic components are
bounded. Thus, with $f_{j,r}=F_{j,r}\circ\chi^{-1}$,
\eqref{eq:Cayley-unitary} gives the bounded cross block
\begin{equation}\label{eq:phase-block}
 U_+N_{f_{j,r}}U_-^{-1}=P M_{\sigma F_{j,r}}|_{\overline{A^2(\D)}}.
\end{equation}
The two Cayley unitaries act on separate spaces, so the coanalytic
projection in the kernel pairing is the full $\bar P$, including constants.
Block multiplication and analytic reproduction now give the right-hand
side of \eqref{eq:compression} with all components dilated.

Let $r\rightarrow1$. The full symbols converge pointwise with a common
uniform bound, so the Toeplitz factors and their products converge
strongly. The disk components converge in $H^2(\D)$ and hence in
$L^2(\D,dA)$. For fixed $z,w$, multiplication by $\sigma$ and by the
fixed Bergman kernels preserves this $L^2(\D,dA)$ convergence. The projected
kernel pairings therefore converge, while the products
$H_{j,r}\overline{G_{j,r}}$ converge in $L^1(\D,dA)$ by Cauchy--Schwarz.
Writing $\mathcal D$ for the analytic compression of
\eqref{eq:halfplane-zero} transported by $U_+$, we obtain
\begin{align}\label{eq:compression}
 (\mathcal D K_w)(z)
 ={}&\frac{\displaystyle
 \sum_{j=1}^NF_j(z)H_j(z)
 +\overline{\sum_{j=1}^NG_j(w)K_j(w)}
 +\sum_{j=1}^NF_j(z)\overline{K_j(w)}}{(1-z\bar w)^2}\notag\\
 &+\sum_{j=1}^N\mathcal C_{\sigma F_j,\sigma K_j}(z,w)
 +\mathcal W_{\sum_{j=1}^NH_j\overline{G_j}}(z,w).
\end{align}
By hypothesis, $\mathcal D=0$. Multiply \eqref{eq:compression}
by $(1-z\bar w)^4$. The first line has finite separation rank,
and \eqref{eq:weighted-cross-kernel} handles the cross terms.
Equation~\eqref{eq:weighted-rigidity} yields
$\sum_{j=1}^NH_j'\overline{G_j'}=0$.
Since $\chi'$ has no zeros, returning to the half-plane gives
\begin{equation}\label{eq:mixed1}
 \sum_{j=1}^N h_j'\overline{g_j'}=0.
\end{equation}
Since the harmonic Bergman projection commutes with pointwise complex
conjugation, \eqref{eq:halfplane-zero} also holds with $u_j$ and $v_j$
replaced by $\overline{u_j}$ and $\overline{v_j}$, respectively.
Applying \eqref{eq:mixed1} to $\overline{u_j}=g_j+\overline{f_j}$ and $\overline{v_j}=k_j+\overline{h_j}$ gives $\sum_{j=1}^N k_j'\overline{f_j'}=0$.
Taking complex conjugates yields
\begin{equation}\label{eq:mixed2}
 \sum_{j=1}^N f_j'\overline{k_j'}=0.
\end{equation}
Uniqueness from the diagonal for holomorphic--antiholomorphic
functions extends these identities to independent variables:
\begin{equation}\label{eq:polarized}
 \sum_{j=1}^Nh_j'(z)\overline{g_j'(w)}=0,\qquad
 \sum_{j=1}^Nf_j'(z)\overline{k_j'(w)}=0.
\end{equation}
This uniqueness follows directly from local power series: all$(z-a)^s \overline{(z-a)}^t(s,t\ge0)$ coefficients of the diagonal identity vanish,
and separate analytic continuation then extends the resulting
identity in the two independent variables.

Fix $\lambda\in\HH$ and put
\[
 \Psi(z,w)=\sum_{j=1}^Nu_j(z)v_j(w),\qquad
 \Psi_0(z,w)=\Psi(z,w)-\Psi(z,\lambda)-\Psi(\lambda,w)+\Psi(\lambda,\lambda).
\]
Integrating \eqref{eq:polarized} separately in the two variables
eliminates the two mixed parts of $\Psi_0$. Thus $\Psi_0$ is a
sum of a purely analytic tensor and a purely coanalytic tensor.
We can choose bounded analytic factors for these remaining tensors,
as follows. Let $X_+$ and $X_-$ be the spans of
$f_j-f_j(\lambda)$ and $\overline{g_j-g_j(\lambda)}$, respectively,
and let $Y_+$ and $Y_-$ be the spans of
$h_j-h_j(\lambda)$ and $\overline{k_j-k_j(\lambda)}$, respectively.
All four spaces consist of functions vanishing at $\lambda$, so
\[
 X_+\cap X_-=\{0\},\qquad Y_+\cap Y_-=\{0\},\qquad
 \Psi_0\in X_+\otimes Y_+ +X_-\otimes Y_-.
\]
Choose a basis of $Y_+\oplus Y_-$ consisting of bases of the
two summands. Finitely many point evaluations span its dual:
otherwise a nonzero function in this space would be annihilated
by every point evaluation. The first-variable coefficients in
this basis are therefore constant linear combinations of finitely
many slices $\Psi_0(z,w_i)$. They are bounded and belong to
$X_+$ or $X_-$, respectively. Choose bases again in the spans
of the nonzero coefficients of each type. Applying the same
evaluation argument in the first variable shows that the
second-variable coefficients are also bounded, and remain in
$Y_+$ or $Y_-$. Each supremum norm is bounded by
$\|\Psi_0\|_\infty$ times a finite sum of the absolute values
of the evaluation coefficients. The resulting boundedness is
therefore global on $\HH$.

We have obtained finite families
$\alpha_\ell,\gamma_\ell,\delta_\ell,e_\ell\in H^\infty(\HH)$ and
$\alpha,\gamma\in h^\infty(\HH)$ such that
\begin{equation}\label{eq:bounded-representation}
 \Psi(z,w)=\alpha(z)+\gamma(w)
 +\sum_\ell \alpha_\ell(z)\gamma_\ell(w)
 +\overline{\sum_\ell \delta_\ell(z)e_\ell(w)}.
\end{equation}
The additive terms can be
taken to be $\alpha(z)=\Psi(z,\lambda)$ and
$\gamma(w)=\Psi(\lambda,w)-\Psi(\lambda,\lambda)$.

Equality of symbol tensors implies equality of the corresponding
finite sums of Toeplitz products. Indeed, if
$\sum_{j=1}^Np_j(z)q_j(w)=0$, choose a basis $e_\ell$ of the span of
the $p_j$ and write $p_j=\sum_\ell c_{\ell j}e_\ell$.
Comparison of basis coefficients gives
$\sum_{j=1}^Nc_{\ell j}q_j=0$, so bilinearity implies
$\sum_{j=1}^NT_{p_j}T_{q_j}=0$. This proves that the tensor-to-operator
map is well defined.

Let $\pi(z)=\Psi(z,z)$. Formula
\eqref{eq:bounded-representation} shows that $\pi(z)$ is
bounded and harmonic. By the block formula \eqref{eq:block},
the analytic diagonal compression of the corresponding operator
is exactly the analytic Toeplitz operator $\TT_{\pi(z)}^{\HH}$.
Thus $\TT_{\pi(z)}^{\HH}=0$, and consequently
\begin{equation}\label{eq:diagonal-zero}
 \pi(z)=0.
\end{equation}
Only injectivity of the analytic Toeplitz map is needed here.
The analytic Cayley unitary identifies this operator with the
disk Toeplitz operator with symbol $\pi\circ\chi$.
Vanishing of all its matrix coefficients makes every mixed
moment of this symbol zero. Polynomials in $z,\bar z$ are dense
in $L^2(\D, dA)$, so a bounded symbol with this property is zero.

By \eqref{eq:bounded-representation} and
\eqref{eq:diagonal-zero}, the original operator equals
\begin{align}\label{eq:semicommutators}
 &\sum_\ell(T_{\alpha_\ell}^{\HH}T_{\gamma_\ell}^{\HH}
                        -T_{\alpha_\ell \gamma_\ell}^{\HH})\notag+\sum_\ell(T_{\bar \delta_\ell}^{\HH}T_{\bar e_\ell}^{\HH}
                        -T_{\overline{\delta_\ell e_\ell}}^{\HH}).
\end{align}
Since $M_\alpha M_\gamma=M_{\alpha\gamma}$ and $L_\alpha L_\gamma=L_{\alpha\gamma}$, the first sum
occupies only the cross block from the coanalytic space to the
analytic space. The second sum occupies the opposite cross block.
Hence each sum is zero.

Consider the first sum, and choose primitives $A_\ell'=\alpha_\ell$
and $B_\ell'=\gamma_\ell$. Because
\[
 L_\gamma\overline{K_w^{\HH}}=\gamma(w)\overline{K_w^{\HH}},
\]
Lemma~\ref{lem:H-cross} gives the cross kernel of an individual
semicommutator as
\[
 \alpha(z)C_\gamma(z,w)+\gamma(w)C_\alpha(z,w)-C_{\alpha\gamma}(z,w).
\]
Expansion and summation give $2H(z,w)/(z-w)^3$, where
\begin{align}\label{eq:H-numerator}
 H(z,w)={}&(z-w)\sum_\ell \alpha_\ell(z)\gamma_\ell(w)\notag-\sum_\ell \alpha_\ell(z)[B_\ell(z)-B_\ell(w)]\notag\\
 &-\sum_\ell \gamma_\ell(w)[A_\ell(z)-A_\ell(w)]\notag+\int_w^z\sum_\ell \alpha_\ell(\xi)\gamma_\ell(\xi)\,d\xi.
\end{align}
Differentiating first in $z$, the undifferentiated terms cancel,
leaving
\[
 H_z(z,w)=\sum_\ell \alpha_\ell'(z)
 [(z-w)\gamma_\ell(w)-B_\ell(z)+B_\ell(w)].
\]
A further derivative in $w$ gives
\begin{equation}\label{eq:H-derivative}
 H_{zw}(z,w)=(z-w)\sum_\ell \alpha_\ell'(z)\gamma_\ell'(w).
\end{equation}
The cross block is zero, so all its kernel columns vanish.
Therefore $H=0$, and
\[
 \sum_\ell \alpha_\ell'(z)\gamma_\ell'(w)=0.
\]
Integrating from $\lambda$ in each variable yields
\[
 \sum_\ell[\alpha_\ell(z)-\alpha_\ell(\lambda)][\gamma_\ell(w)-\gamma_\ell(\lambda)]=0.
\]
The opposite cross block is treated in the same way. Thus
$\Psi_0=0$, and $\Psi(z,w)=\alpha_0(z)+\gamma_0(w)$.
Equation~\eqref{eq:diagonal-zero} gives
$\gamma_0=-\alpha_0$. Taking $\phi(z)=\Psi(z,\lambda)$ gives a function
in $h^\infty(\HH)$ with $\phi(\lambda)=0$ and proves
\eqref{eq:additive}.

Conversely, \eqref{eq:additive} and the well-defined
tensor-to-operator map imply
\[
 \sum_{j=1}^NT_{u_j}^{\HH}T_{v_j}^{\HH}
 =T_\phi^{\HH}T_1^{\HH}-T_1^{\HH}T_\phi^{\HH}=0.\qedhere
\]
\end{proof}

\subsection{Boundary limits}
\label{subsec:compact-limits}

We apply affine rescaling to the Toeplitz sum $D$ in \eqref{eq:S},
using the tensor $\Phi$ and the analytic decompositions
\eqref{eq:components}. No individual analytic component is assumed
to be bounded. Composition with the affine maps defined below preserves harmonicity
on the rescaled domains. As these domains exhaust the half-plane,
we can study the boundary behavior of the symbols through their
rescaled limits.

\begin{lemma}\label{lem:limits}
For any sequence $a_n=r_n\zeta_n$ with $r_n\to1$ and
$|\zeta_n|=1$, put
\begin{equation}\label{eq:scaling}
 \delta_n=1-r_n,\qquad L_n(z)=\zeta_n(1-\delta_nz),\qquad
 \Omega_n=L_n^{-1}(\D)=\{z:2\operatorname{Re}z>\delta_n|z|^2\}.
\end{equation}
For each subsequence along which all the rescaled symbols converge,
retain the index $n$ and write
\begin{equation}\label{eq:symbol-limits}
 \varphi_j\circ L_n\longrightarrow u_j,\qquad
 \psi_j\circ L_n\longrightarrow v_j
 \quad\text{in }\Loc(\HH),
\end{equation}
and set $D_\infty=\sum_{j=1}^NT_{u_j}^{\HH}T_{v_j}^{\HH}$.
Then $D$ is compact if and only if
\begin{equation}\label{eq:limit-criterion}
 D_\infty=0
\end{equation}
for every sequence of centers tending to the boundary and every
simultaneous symbol limit along a subsequence of that sequence.
\end{lemma}
\begin{proof}
Interior estimates for bounded harmonic functions show that every
boundary sequence has a subsequence satisfying
\eqref{eq:symbol-limits}. The limiting symbols are bounded and
harmonic. Define the isometry
\begin{equation}\label{eq:scale-isometry}
 U_n:L^2(\D, dA)\longrightarrow L^2(\HH, dA),\qquad
 U_nf=\mathbf1_{\Omega_n}\delta_n(f\circ L_n),
\end{equation}
and put $Q_n=U_nQU_n^*$. Since $L_n$ maps $\Omega_n$ onto $\D$ and $|L_n'|=\delta_n$,
the change of variables $w=L_n(z)$ gives
\[
 \|U_n f\|_{L^2(\HH,dA)}^2
 =\delta_n^2\int_{\Omega_n}|f(L_n(z))|^2\,dA(z)
 =\int_{\D}|f(w)|^2\,dA(w)
 =\|f\|_{L^2(\D,dA)}^2.
\]
For every $F\in L^2(\HH,dA)$, the rescaling
argument of \cite[proof of Lemma~3.4]{CHLZ} for these domains and isometries gives
\begin{equation}\label{eq:projection-strong}
 Q_n\longrightarrow Q_{\HH}\quad\text{strongly on }L^2(\HH,dA),
 \qquad U_n^*F\rightharpoonup0\quad\text{in }L^2(\D,dA),
\end{equation}
as well as the normalized-kernel convergence
\begin{equation}\label{eq:normalized-kernels}
 U_nr_{a_n}\longrightarrow r_1^{\HH}
 \quad\text{in }L^2(\HH).
\end{equation}
Here $r_1^{\HH}$ is the normalized harmonic reproducing kernel
at $1$ on the half-plane. These facts depend only on the
rescaled spaces and kernels. We apply them to the present finite sum.

Extend every rescaled symbol by zero outside $\Omega_n$.
Uniform boundedness and dominated convergence give strong
convergence of its multiplication operator. Finite products of
uniformly bounded strongly convergent operators also converge
strongly. Consequently,
\begin{equation}\label{eq:operator-strong}
 U_n(DQ)U_n^*\longrightarrow D_\infty Q_{\HH}
 \qquad\text{strongly}.
\end{equation}
If $D$ is compact, the weak convergence in
\eqref{eq:projection-strong} makes the left side converge strongly
to zero, so $D_\infty=0$.

Conversely, suppose that all boundary limits vanish. Along each
subsequence with simultaneous symbol limits,
\eqref{eq:normalized-kernels}, \eqref{eq:operator-strong}, and
the uniform operator-norm bounds give
\[
 \|Dr_{a_n}\|_2
 =\|U_n(DQ)U_n^*U_nr_{a_n}\|_2\longrightarrow0.
\]
Subsequence extraction yields $\|Dr_a\|_2\to0$ as
$|a|\to1$. For fixed $R$, points of
$D_\beta(a,R)$ approach the boundary uniformly as
$|a|\to1$. Thus
$|\langle Dr_b,r_a\rangle|\le\|Dr_b\|_2$ and
\eqref{eq:localization} prove compactness.

For later use, if $b_n\in D_\beta(a_n,R)$ and
$w_n=L_n^{-1}(b_n)$, then, with $\rho=\tanh R<1$,
\[
 |w_n-1|<\rho|1+(1-\delta_n)w_n|.
\]
The triangle inequality and expansion of the squared inequality give
\[
 |w_n|\le\frac{1+\rho}{1-\rho},\qquad
 \operatorname{Re}w_n>
 \frac{1-\rho^2}{2(1+\rho^2)}.
\]
Hence these rescaled points lie in a fixed compact subset of $\HH$.
\end{proof}

The following finite-jet consequence of the classical Wronskian
theorem \cite{BD} will be used to bound the derivative order.

\begin{lemma}
\label{lem:wronskian}
Let $\Omega\subset\C$ be a domain, let $\alpha_j,\gamma_j$ be
holomorphic on $\Omega$, and set
$\mathcal{R}(z,w)=\sum_{j=1}^{N}\alpha_j(z)\gamma_j(w)$. If
\begin{equation}\label{eq:diagonal-jets}
 (\partial_z^i\mathcal{R})(z,z)=0,\qquad i=0,\ldots,N-1,
\end{equation}
then $\mathcal{R}=0$.
\end{lemma}
\begin{proof}
Choose a basis $A_1,\ldots,A_r$ of the constant-coefficient
linear span of the $\alpha_j$, with $r\le N$. Then
$\mathcal{R}=\sum_{\ell=1}^{r}A_\ell(z)B_\ell(w)$.
If $r=0$, the assertion is immediate. Otherwise the first $r$
identities in \eqref{eq:diagonal-jets} give
\begin{equation}\label{eq:wronskian-system}
 \begin{pmatrix}
 A_1&\cdots&A_r\\ A_1'&\cdots&A_r'\\
 \vdots&&\vdots\\ A_1^{(r-1)}&\cdots&A_r^{(r-1)}
 \end{pmatrix}
 \begin{pmatrix}B_1\\\vdots\\B_r\end{pmatrix}=0.
\end{equation}
By the classical analytic Wronskian theorem
\cite[Theorem~1]{BD}, the Wronskian of $A_1,\ldots,A_r$ is
not identically zero. Thus the matrix in
\eqref{eq:wronskian-system} is invertible on a nonempty open
set. Each $B_\ell$ is zero there and hence, by the identity
theorem, throughout $\Omega$.
\end{proof}

\begin{proof}[Proof of Theorem~\ref{thm:compact}]
\textit{Equivalence of \textup{(i)}--\textup{(iii)}.}
Suppose first that $D$ is compact. By Lemma~\ref{lem:limits}
and Proposition~\ref{prop:halfplane-zero}, the symbol tensor
of each simultaneous boundary limit has the form
\begin{equation}\label{eq:limit-additive}
 \Psi(z,w):=\sum_{j=1}^Nu_j(z)v_j(w)=\phi(z)-\phi(w).
\end{equation}
In particular, $\Psi(z,z)=0$ and
\begin{equation}\label{eq:limit-centered}
 \Psi(z,w)-\Psi(z,1)-\Psi(1,w)+\Psi(1,1)=0.
\end{equation}
If $\Phi(a,a)\not\to0$, choose a boundary sequence contradicting
this limit and then extract a simultaneous symbol limit.
Since $L_n(1)=a_n$, we obtain
$\Phi(a_n,a_n)\to\Psi(1,1)=0$, a contradiction.

Fix $R>0$. If $\omega_R(a)\not\to0$, we can choose
$a_n\to\partial\D$ and
$z_n,w_n\in D_\beta(a_n,R)$ such that the modulus of the
centered expression has a positive lower bound. The estimates
in the proof of Lemma~\ref{lem:limits} place both rescaled
points in a fixed compact subset of $\HH$. Pass to convergent
subsequences and use simultaneous local uniform convergence of
the symbols. This contradicts \eqref{eq:limit-centered}.
Thus condition \textup{(iii)} holds, and \textup{(iii)} plainly
implies \textup{(ii)}.

Finally, suppose that \textup{(ii)} holds for some $R>0$.
Take any simultaneous boundary limit $\Psi$. For fixed
$z\in\HH$, we have $|L_n(z)|\to1$, so
\begin{equation}\label{eq:limit-diagonal}
 \Psi(z,z)=\lim_n \Phi(L_n(z),L_n(z))=0.
\end{equation}
Choose a sufficiently small open neighborhood $V\Subset\HH$
of $1$ such that $L_n(V)\subset D_\beta(a_n,R)$ for all
large $n$. Such a neighborhood exists because
\[
 \left|\frac{L_n(z)-a_n}{1-\bar a_nL_n(z)}\right|
 =\frac{|z-1|}{|1+(1-\delta_n)z|}
\]
converges locally uniformly near $1$. The centered condition
gives \eqref{eq:limit-centered} on $V\times V$.
First fix $w\in V$ and apply unique continuation for harmonic
functions in the first variable. Then fix any $z\in\HH$ and
repeat the argument in the second variable. The identity follows
on all of $\HH^2$.

Consequently, $\Psi(z,w)=\alpha(z)+\beta(w)$, where
$\alpha,\beta$ can be chosen as bounded harmonic slices.
Equation~\eqref{eq:limit-diagonal} gives $\beta=-\alpha$.
Proposition~\ref{prop:halfplane-zero} yields $D_\infty=0$,
and Lemma~\ref{lem:limits} implies that $D$ is compact.

\textit{Equivalence of \textup{(i)} and \textup{(iv)}.}
Use the rescaling from Lemma~\ref{lem:limits}, and write
$\lambda_n=-\delta_n\zeta_n$, so that $L_n'=\lambda_n$.
For $i\ge1$,
\begin{align}
 \partial_z^i(\varphi_j\circ L_n)
 &=\lambda_n^iF_j^{(i)}\circ L_n,\label{eq:scale-hol}\\
 \partial_{\bar z}^i(\varphi_j\circ L_n)
 &=\bar\lambda_n^i\overline{G_j^{(i)}\circ L_n}.
 \label{eq:scale-anti}
\end{align}
The same formulas hold for $\psi_j$. Moreover,
\begin{equation}\label{eq:scale-weight}
 1-|L_n(z)|^2
 =\delta_n(2\operatorname{Re}z-\delta_n|z|^2).
\end{equation}
On each fixed compact subset of the half-plane,
$1-|L_n(z)|^2$ is therefore comparable to $\delta_n$.
Together with convergence in $\Loc$, this turns the weighted
boundary limits of the original symbols into exact derivative
identities for the limiting symbols.

Suppose first that $D$ is compact. By
Proposition~\ref{prop:halfplane-zero} and
Lemma~\ref{lem:limits}, each limiting tensor is an additive
difference. Consequently,
\begin{equation}\label{eq:necessary-limit-diagonal}
 \sum_{j=1}^Nu_jv_j=0,
\end{equation}
\begin{equation}\label{eq:necessary-limit-mixed}
 \sum_{j=1}^N(\partial_zu_j)(z)(\partial_{\bar w}v_j)(w)=0,
\end{equation}
and, for every $i\ge1$,
\begin{align}
 \sum_{j=1}^N(\partial_z^iu_j)(z)(\partial_{ w}v_j)(w)&=0,
 \label{eq:necessary-limit-hol}\\
 \sum_{j=1}^N(\partial_{\bar z}^iu_j)(z)(\partial_{\bar w}v_j)(w)&=0.
 \label{eq:necessary-limit-anti}
\end{align}
If any of the limits in
\eqref{eq:d-vanish}--\eqref{eq:anti-vanish} fails, choose a
boundary sequence along which its modulus has a positive lower
bound. Rescale about that sequence and extract simultaneous
symbol limits. At $z=1$,
$(1-|a_n|^2)/\delta_n=2-\delta_n\to2$ and
$|\lambda_n|=\delta_n$. By
\eqref{eq:scale-hol}--\eqref{eq:scale-weight}, the modulus
of the corresponding limiting derivative expression retains a
positive lower bound. This contradicts
\eqref{eq:necessary-limit-diagonal}--\eqref{eq:necessary-limit-anti}.
Hence all the stated boundary conditions hold.

Conversely, assume
\eqref{eq:d-vanish}--\eqref{eq:anti-vanish}. Take any
simultaneous boundary limits $u_j,v_j$ and set
$\Psi(z,w)=\sum_{j=1}^Nu_j(z)v_j(w)$.
Condition~\eqref{eq:d-vanish} gives
\begin{equation}\label{eq:limit-d}
 \sum_{j=1}^Nu_j(z)v_j(z)=0.
\end{equation}
Condition~\eqref{eq:mixed-vanish} and the scaling formulas give
\begin{equation}\label{eq:limit-mixed-diag}
 \sum_{j=1}^N(\partial_zu_j)(z)(\partial_{\bar z}v_j)(z)=0.
\end{equation}
Diagonal uniqueness \cite[Lemma~3.1]{CKL} therefore yields
\begin{equation}\label{eq:mixed-tensor1}
 \partial_z\partial_{\bar w}\Psi(z,w)
 =\sum_{j=1}^N(\partial_zu_j)(z)(\partial_{\bar w}v_j)(w)=0.
\end{equation}
Apply $\partial\bar\partial$ to the exact limit identity
\eqref{eq:limit-d}. Harmonicity and
\eqref{eq:limit-mixed-diag} give
$\sum_{j=1}^N(u_j)_{\bar z}(z)(v_j)_z(z)=0$.
Applying the same diagonal uniqueness to this identity gives
\begin{equation}\label{eq:mixed-tensor2}
 \partial_{\bar z}\partial_w\Psi(z,w)=0.
\end{equation}
Consider the holomorphic function
\[
 F_+(z,w)= \partial_z\partial_w\Psi(z,w)
 =\sum_{j=1}^N(\partial_zu_j)(z)(\partial_wv_j)(w).
\]
Its separation rank is at most $N$. Condition
\eqref{eq:hol-vanish} and the scaling formulas give, for
$i=1,\ldots,N$,
\[
 \sum_{j=1}^N(\partial_z^iu_j)(z)(\partial_zv_j)(z)=0,
 \quad\text{equivalently}\quad
 (\partial_z^{i-1}F_+)(z,z)=0.
\]
Lemma~\ref{lem:wronskian} implies
\begin{equation}\label{eq:pure-tensor1}
 \partial_z\partial_w\Psi=0.
\end{equation}
Take complex conjugates in \eqref{eq:anti-vanish}, and use
\eqref{eq:scale-anti}. Applying the same lemma to the
holomorphic function
$\overline{\partial_{\bar z}\partial_{\bar w}\Psi(z,w)}$ gives
\begin{equation}\label{eq:pure-tensor2}
 \partial_{\bar z}\partial_{\bar w}\Psi=0.
\end{equation}

Equations~\eqref{eq:mixed-tensor1}--\eqref{eq:pure-tensor2}
say that all four types of cross-variable mixed Wirtinger
derivatives of $\Psi$ vanish. Thus
$\Psi(z,w)=\alpha(z)+\beta(w)$. To verify this directly,
fix $o\in\HH$ and set
\[
 K(z,w)=\Psi(z,w)-\Psi(z,o)-\Psi(o,w)+\Psi(o,o).
\]
Both first derivatives of $\Psi_z$ and $\Psi_{\bar z}$ in
the second variable vanish, so these functions are independent
of that variable. Hence $\partial_zK=\partial_{\bar z}K=0$.
Since $K(o,w)=0$, it follows that $K=0$.
The functions $\alpha,\beta$ can therefore be chosen as bounded
harmonic slices. Equation~\eqref{eq:limit-d} gives
$\beta=-\alpha$. Proposition~\ref{prop:halfplane-zero}
and Lemma~\ref{lem:limits} show that $D$ is compact.
\end{proof}

\begin{remark}\label{rem:product-types}
The reflection $U$ preserves area and the harmonic subspace, so $UQ=QU$
and $UT_{\gamma}U=T_{U\gamma}$. Since $h_{\alpha}=T_{\alpha}U$,
\[
 h_{\alpha}h_{\gamma}=T_{\alpha}T_{U\gamma},\qquad (T_{\alpha}h_{\gamma})U=T_{\alpha}T_{\gamma},\qquad (h_{\alpha}T_{\gamma})U=T_{\alpha}T_{U\gamma}.
\]
Thus any sum of $TT$ and $hh$ terms is itself a Toeplitz product sum. A sum of $Th$ and $hT$ terms becomes a Toeplitz
product sum upon right multiplication by $U$, which preserves rank,
compactness and essential norm. Since $U^2=I$,
both substitutions are reversible. A sum containing terms from both groups
has the form $D_0+D_1U$, with $D_0,D_1$ as in \eqref{eq:S}. The two main
theorems give criteria for $D_0$ and $D_1U$ separately.
\end{remark}

\begin{remark}\label{rem:compact-even-order}
For every even $N\ge2$, take $p(z)=(1-z)^\mathrm{i}$ with the principal
logarithm and, for $0\le j<N$, put
\[
 \varphi_j=\frac{(-1)^j\binom{N-1}{j}}{(N-j)(j+1)}p^{N-j},
 \qquad \psi_j=p^{j+1}.
\]
These symbols are bounded and analytic since $|\arg(1-z)|<\pi/2$.
Pairing $j$ with $N-1-j$ gives $\sum_j\varphi_j\psi_j=0$, while
the binomial theorem gives
\[
 \sum_j\varphi_j'(z)\psi_j'(w)
 =p'(z)p'(w)\bigl(p(z)-p(w)\bigr)^{N-1}.
\]
Differentiate $i-1$ times in $z$ before putting $w=z$.
The result vanishes for $1\le i<N$ and equals $(N-1)!(p'(z))^{N+1}$
for $i=N$. Along $z=r\rightarrow1$ its modulus times $(1-r^2)^{N+1}$
tends to $(N-1)!\,2^{N+1}$. All mixed and coanalytic conditions vanish
identically. Thus omitting order $N$ in Theorem~\ref{thm:compact}\textup{\,(iv)}
would admit this noncompact sum. The case $N=2$ recovers, up to a scalar,
\cite[Example~4.10(i)]{CHLZ}.
\end{remark}

\section{Consequences and examples}\label{sec:examples}

\begin{corollary}\label{cor:general-polynomial}
Let $D$ be the sum in \eqref{eq:S}.
\begin{enumerate}
\item If all $\varphi_j,\psi_j$ are harmonic polynomials, then $D$
has finite rank if and only if
\[
 \sum_{j=1}^{N}\varphi_j(z)\psi_j(z)=0,\qquad
 \sum_{j=1}^{N}(\partial_z\varphi_j)(z)
                    (\partial_{\bar z}\psi_j)(z)=0
 \qquad(z\in\D).
\]
\item For arbitrary bounded harmonic symbols, the tensor $\Phi$ in
\eqref{eq:Phi} characterizes the zero operator:
\[
 D=0 \quad\Longleftrightarrow\quad
 \Phi(z,w)=\phi(z)-\phi(w)\quad(z,w\in\D)
 \quad\text{for some }\phi\in\hinfty .
\]
If the two symbol identities in \textup{(i)} hold, then $D=0$ also
holds if and only if both combined numerators $\mathcal{R}_+$ and $\mathcal{R}_-$ vanish
identically.
\end{enumerate}
\end{corollary}
\begin{proof}
Diagonal uniqueness and integration identify the mixed derivative
condition with the centered two-variable identity in
Theorem~\ref{thm:general}. For polynomial components, the two
numerators are polynomials and their quotients by $(z-w)^3$ are
holomorphic by \eqref{eq:Rintegral}. Hence both numerators are
divisible by $(z-w)^3$. The rational bases can be chosen as
$\{1\}$, so \eqref{eq:rd3} is automatic. The assertion about
$\mathcal{R}_+$ and $\mathcal{R}_-$ follows from the two cross kernels in
Lemma~\ref{lem:cancellation}.

Suppose $D=0$. Proposition~\ref{prop:input} gives $\Phi(z,z)=0$
and the two centered identities, hence $\Phi_{z\bar w}=\Phi_{\bar z w}=0$.
Also $\mathcal{R}_+=\mathcal{R}_-=0$. Differentiating \eqref{eq:Rdef} gives
$\partial_z\partial_w\mathscr R_{f,h}(z,w)=(z-w)f'(z)h'(w)$.
Applied to $\mathcal{R}_+$ at $(z,w)$ and to $\mathcal{R}_-$ at $(\bar z,\bar w)$, this gives
$\Phi_{zw}=\Phi_{\bar z\bar w}=0$ off the diagonal and hence everywhere.
Thus $\Phi(z,w)-\Phi(z,0)$ is independent of $z$ and equals
$\Phi(0,w)-\Phi(0,0)$. The diagonal identity now gives the claimed
representation with $\phi(z)=\Phi(z,0)\in\hinfty$.
Conversely, bilinearity on symbol tensors, as in the proof of
Proposition~\ref{prop:halfplane-zero}, gives $D=T_\phi T_1-T_1T_\phi=0$.
\end{proof}

Remark~\ref{rem:product-types} reduces each of the following commutator problems to a difference of Toeplitz products. 
The ordinary Toeplitz commutator is classified in
\cite[Theorems~1.1 and~1.2, Corollary~4.1]{CHLZ}. The distinction is already
visible for $\varphi(z)=z$ and $\psi=1$:
\[
 [T_z,T_1]=0,\qquad [T_z,h_1]U=[h_z,h_1]=T_{z-\bar z}.
\]
The last operator is not compact by Theorem~\ref{thm:compact}, since
its diagonal symbol $z-\bar z$ tends to $2i$ along $z=\mathrm{i}r$, $r\rightarrow1$.
Indeed every ordinary Toeplitz commutator with harmonic symbols annihilates
$1$, since $[T_f,T_g]1=Q(fg)-Q(gf)$, whereas $T_{z-\bar z}1=z-\bar z$.

Here we consider pairs of bounded harmonic symbols. For
$\varphi,\psi\in h^\infty(\D)$, reflection gives
\[
 [T_\varphi,h_\psi]U
 =[T_\varphi,T_\psi]+T_\psi T_{\varphi-U\varphi}.
\]
The mixed commutator has finite rank if and only if it is zero.
These conditions hold precisely when $\psi=0$, or when
$\varphi=U\varphi$ and $1,\varphi,\psi$ are linearly
dependent over $\C$. Indeed, finite rank and
Proposition~\ref{prop:input} give
$\psi(\varphi-U\varphi)=0$. If $\psi\ne0$, unique
continuation gives $\varphi=U\varphi$, and
\cite[Theorem~1.1]{CHLZ} gives the stated dependence, since a
reflection-invariant analytic or coanalytic function is constant.
Conversely, the mixed commutator vanishes if $\psi=0$ or
$\varphi$ is constant. In the remaining case,
$\psi=\lambda\varphi+\mu$ and $UT_\varphi=T_\varphi U$,
which again give a zero commutator.

Moreover, $[T_\varphi,h_\psi]$ is compact if and only if
$[T_\varphi,T_\psi]$ is compact and
\[
 \psi(z)\bigl[\varphi(z)-\varphi(\bar z)\bigr]
 \longrightarrow0\qquad(|z|\to1).
\]
Necessity of this boundary condition follows from
Theorem~\ref{thm:compact}\textup{\,(ii)}. Under this condition,
every simultaneous half-plane profile $(\Psi,V)$ of
$(\psi,\varphi-U\varphi)$ satisfies $\Psi V=0$.
Unique continuation makes one profile identically zero, so
Lemma~\ref{lem:limits} gives compactness of
$T_\psi T_{\varphi-U\varphi}$. The displayed decomposition
and the unitarity of $U$ prove the equivalence. The remaining
Toeplitz commutator is compact precisely under the three weighted
derivative conditions in \cite[Theorem~1.2]{CHLZ}.
If both symbols extend continuously to $\overline\D$, the Toeplitz
commutator is automatically compact. In this case the remaining
boundary condition also follows from the quotient-algebra description
in \cite[Theorem~4.4]{GZ}, using $h_\psi=T_\psi U$.

The reflected self-commutator from \cite[Example~4.9]{CHLZ} gives
\[
 [T_{z+\bar z},h_z]=[T_{\bar z},T_z]U,
 \qquad \lambda_n=\frac1{(n+1)(n+2)}\quad(n\ge1).
\]
Since $Uz^n=\overline{z}^n$ and $U\overline{z}^n=z^n$, the cited diagonal
action becomes
\[
 \begin{aligned}
 [T_{z+\bar z},h_z]1 &= 0,
 & [T_{z+\bar z},h_z]z^n &= -\lambda_n\bar z^n,\\
 [T_{z+\bar z},h_z]\bar z^n &= \lambda_n z^n,
 & [T_{z+\bar z},h_z]^2z^n &= -\lambda_n^2z^n.
 \end{aligned}
\]
The orthogonal two-dimensional blocks have singular values
$\lambda_n\to0$, so $[T_{z+\bar z},h_z]$ is compact of infinite rank and $[T_{z+\bar z},h_z]^2\ne0$.

Chen, He, and Hu \cite{CHH} characterized commutativity of small
Hankel operators with bounded harmonic symbols.
We now classify the bounded harmonic symbol pairs that give positive
finite rank on the full harmonic Bergman space.

\begin{corollary}\label{cor:hh-finite}
Let $\varphi,\psi\in h^\infty(\D)$. Then $[h_\varphi,h_\psi]$ has
nonzero finite rank if and only if there is a nonconstant analytic
function $u$ such that $u,u^{-1}\in H^\infty(\D)$, both have rational
primitives, and for some $a,b,c,d\in\C$,
\begin{equation}\label{hh:reciprocal-pair}
 \varphi=a u+b\,U(u^{-1}),\qquad
 \psi=c u+d\,U(u^{-1}),\qquad ad-bc\ne0.
\end{equation}
In this case
\begin{equation}\label{hh:even-rank}
 \rank[h_\varphi,h_\psi]=2\rank (T_u,T_{u^{-1}}].
\end{equation}
Every positive even rank occurs.
\end{corollary}

The evenness of the rank also follows from the theorem of Chen, Koo,
and Lee on finite-rank complex skew-symmetric operators
\cite{ChenKooLeeSkew}. Define the conjugation $K$ on the harmonic
Bergman space by $(Kf)(z)=\overline{f(\bar z)}$. Since
$Kh_\alpha^*K=h_\alpha$, we have
\[
K[h_\varphi,h_\psi]^*K=-[h_\varphi,h_\psi].
\]
We now prove the symbol classification and the rank formula, and show
that every positive even rank occurs.

\begin{proof}
Suppose $[h_\varphi,h_\psi]\ne0$ has finite rank, and write
$\varphi=\varphi(0)+f+g(\bar z)$,
$\psi=\psi(0)+p+q(\bar z)$ with centered components.
Applying Proposition~\ref{prop:input} and Remark~\ref{rem:one-mixed} to
$[h_\varphi,h_\psi]=T_\varphi T_{U\psi}-T_\psi T_{U\varphi}$ gives
\begin{equation}\label{hh:necessary}
\begin{aligned}
 \varphi(z)\psi(\bar z)-\psi(z)\varphi(\bar z)&=0,\\
 f(z)p(w)-p(z)f(w)&=0,\qquad
 g(z)q(w)-q(z)g(w)=0.
\end{aligned}
\end{equation}
Since $[h_\varphi,h_\psi]\ne0$, $V=\operatorname{span}_{\C}\{\varphi,\psi\}$ is two-dimensional.
If $1\in V$, the first identity, in a basis $(\alpha_0,1)$, gives
$\alpha_0=U\alpha_0$. But
\[
 [h_{\alpha_0},h_1]=T_{\alpha_0-U\alpha_0}=0,
\]
so bilinearity makes $[h_\varphi,h_\psi]=0$, a contradiction.
Hence $1\notin V$. The two centered-component maps on $V$ have
rank at most one and jointly are injective; both therefore have
rank one and distinct kernels. Choosing their kernel vectors gives
a basis $(u,Uv)$ with nonconstant $u,v\in H^\infty$.
Their boundedness follows because these basis vectors belong to $V$.
The first identity says that $uv$ equals its reflection, hence is
a nonzero constant. Rescaling gives $v=u^{-1}$ and
\eqref{hh:reciprocal-pair}; the coefficient determinant is nonzero
because both symbol pairs are bases of $V$.

Direct multiplication gives
\begin{equation}\label{hh:blocks}
 [h_\varphi,h_\psi]=(ad-bc)\bigl((T_u,T_v]-U(T_v,T_u]U\bigr).
\end{equation}
Equations~\eqref{eq:Rintegral} and \eqref{eq:semicommutator-cross-kernel} give
$\mathscr K_{v,u}(z,w)=\mathscr K_{u,v}(w,z)$.
Transposition preserves separation rank. The opposite cross blocks
therefore have equal ranks and orthogonal ranges, proving
\eqref{hh:even-rank} and reducing finite rank of $[h_\varphi,h_\psi]$ to that of
$(T_u,T_v]$. Both semicommutators are then finite rank, as is
their difference $[T_u,T_v]$.
By \cite[Theorem~1.1]{CHLZ}, both primitives are rational unless
$1,u,u^{-1}$ are dependent; multiplying such a relation by $u$
would make its nonconstant range satisfy a nonzero quadratic
polynomial, contradicting the open mapping theorem.
Conversely, suppose that $u$ and $v=u^{-1}$ have rational primitives.
The finite-rank small Hankel theorem \cite[Theorem~4.2]{GZsmall2001},
in the disk formulation of \cite[Theorem~2.2 and the paragraph following it]{CHLZ},
shows that $PM_{\alpha}|_{\overline{\A}}$ has finite rank for $\alpha=u,v$.
Each harmonic Toeplitz operator $T_{\alpha}$ is upper triangular in the
analytic--coanalytic decomposition. Since $uv=1$, the diagonal blocks of
$T_uT_v-I$ vanish, and each term in its remaining block contains one
of these finite-rank restrictions. The same holds for $T_vT_u-I$.
Equation~\eqref{hh:blocks} gives the required finite rank.
Moreover, $\mathscr K_{u,v}(z,z)=z^2(u'(z)/u(z))^2/3\not\equiv0$,
so this branch is nonzero.

For every positive even rank, fix $0<a<1$ and an integer $m\ge2$, and set
\begin{equation}\label{hh:even-family}
\begin{gathered}
 A_m(z)=(1-az)^m,\qquad B_m(z)=(1-az)^{-m}.
\end{gathered}
\end{equation}
Formula~\eqref{eq:semicommutator-cross-kernel} yields
\begin{equation}\label{hh:even-kernel}
 \mathscr K_{A_m,B_m}(z,w)
 =\frac{2a^2m}{m^2-1}zw
  \sum_{k=0}^{m-2}(k+1)(m-1-k)
       \frac{(1-az)^k}{(1-aw)^{k+2}}.
\end{equation}
For $m=2$ this kernel is $4a^2zw/[3(1-aw)^2]$, so the cross block
has rank one and the small Hankel commutator has rank two by \eqref{hh:even-rank}.
In general the factor families $z(1-az)^k$ and $w(1-aw)^{-k-2}$,
$0\le k\le m-2$, are separately linearly independent, so
$\rank (T_{A_m},T_{B_m}]=m-1$ and $\rank[h_{A_m},h_{UB_m}]=2(m-1)$.
Real Taylor coefficients give $[h_{A_m},h_{UB_m}]=(T_{A_m},T_{B_m}]-(T_{A_m},T_{B_m}]^*\ne0$,
hence
\[
 -[h_{A_m},h_{UB_m}]^2
 =[h_{A_m},h_{UB_m}]^*[h_{A_m},h_{UB_m}]\ne0.\qedhere
\]
\end{proof}

\smallskip\noindent\textit{Compactness.}
Apply Theorem~\ref{thm:compact} with $N=2$ and the symbol tensor
\[
 \Phi_h(z,w)=\varphi(z)\psi(\bar w)-\psi(z)\varphi(\bar w).
\]
Its derivative criterion says that $[h_\varphi,h_\psi]$ is compact if and
only if the following four quantities tend to zero as $|z|\to1$:
\[
\begin{gathered}
 \Phi_h(z,z),\qquad (1-|z|^2)^2(\partial_z\partial_{\bar w}\Phi_h)(z,z),\\
 (1-|z|^2)^2(\partial_z\partial_w\Phi_h)(z,z),\qquad
 (1-|z|^2)^3(\partial_z^2\partial_w\Phi_h)(z,z).
\end{gathered}
\]
Here the derivatives are taken before setting $w=z$.
To see why four conditions suffice, put $F(z)=(\partial_z\partial_w\Phi_h)(z,z)$
and $G(z)=(\partial_z^2\partial_w\Phi_h)(z,z)$. By
$\Phi_h(z,w)=-\Phi_h(\bar w,\bar z)$, the two remaining coanalytic
expressions are $-F(\bar z)$ and $G(\bar z)-F'(\bar z)$. Since $F$ is
analytic, Cauchy's estimate on $|w-z|=(1-|z|)/2$ gives
$(1-|z|^2)^3 F'(z)\to0$ whenever $(1-|z|^2)^2F(z)\to0$.

For symbols continuous on $\overline\D$, compactness is already
equivalent to $\Phi_h(z,z)\to0$. This special case follows from
\cite[Theorem~4.4]{GZ}, since the image of $[h_\varphi,h_\psi]$ in
the quotient by compact operators is the boundary function
$\varphi(U\psi)-\psi(U\varphi)$.

A useful reciprocal-symbol case also follows from existing results:
\cite[Corollary~4.5(i)]{CHLZ} and Fulsche's compact semicommutator
result \cite[Lemma~13 in the arXiv version]{Fulsche} give, for $u,v=u^{-1}\in H^\infty(\D)$,
\[
 [h_u,h_{Uv}]\text{ is compact}
 \quad\Longleftrightarrow\quad (1-|z|^2)|u'(z)|\longrightarrow0
 \quad (|z|\to1).
\]
Indeed, \eqref{hh:blocks} has two opposite cross blocks, so compactness
of the small Hankel commutator makes both $(T_u,T_v]$ and $(T_v,T_u]$
compact. Their difference $[T_u,T_v]$ is therefore compact. Since
\[
 u'v''-u''v'=2(u'/u)^3,
\]
\cite[Corollary~4.5(i)]{CHLZ} and the boundedness of $u$ and $u^{-1}$
give the little Bloch condition on the right. Conversely, integration
of $u'$ along hyperbolic geodesics shows that this condition makes the
oscillation of $u$ on every fixed-radius hyperbolic disk tend to zero.
Thus $u\in\mathrm{VO}\subset\mathrm{VMO}_{\partial}$, and Fulsche's
result makes both semicommutators compact. Equation~\eqref{hh:blocks}
then gives the asserted compactness.

\section*{Data availability}
No data were used for the research described in this article.
\section*{Conflict of interest}
The authors declare that they have no conflict of interest.
\section*{Declaration of generative AI and AI-assisted technologies in the manuscript preparation process}
OpenAI's ChatGPT was used to assist with mathematical exposition, to check
illustrative examples, and to revise the English text and LaTeX formatting. The authors independently verified the mathematical arguments, calculations, and references and take full responsibility for the final manuscript.

\end{document}